\documentclass[12pt]{article}
\usepackage{amsmath}
\usepackage{amsthm}
\usepackage{amsfonts}
\usepackage{amssymb}
\newcommand{\be}{\begin{equation}}
\newcommand{\ee}{\end{equation}}
\newcommand{\bea}{\begin{eqnarray}}
\newcommand{\eea}{\end{eqnarray}}
\newcommand{\bean}{\begin{eqnarray*}}
\newcommand{\eean}{\end{eqnarray*}}
\newcommand{\brray}{\begin{array}}
\newcommand{\erray}{\end{array}}
\newcommand{\ben}{\begin{equation}{nonumber}}
\newcommand{\een}{\end{equation}{nonumber}}

\newtheorem{dfn}{Definition}[section]
\newtheorem{thm}[dfn]{Theorem}
\newtheorem{lema}[dfn]{Lemma}
\newtheorem{pro}[dfn]{Proposition}
\newtheorem{coro}[dfn]{Corollary}
\newtheorem{xmpl}[dfn]{Example}
\newtheorem{rmrk}[dfn]{Remark}

\newcommand{\bdfn}{\begin{dfn}}
\newcommand{\bthm}{\begin{thm}}
\newcommand{\blema}{\begin{lema}}
\newcommand{\bpro}{\begin{pro}}
\newcommand{\bcoro}{\begin{coro}}
\newcommand{\bxmpl}{\begin{xmpl}}
\newcommand{\brmrk}{\begin{rmrk}}

\newcommand{\edfn}{\end{dfn}}
\newcommand{\ethm}{\end{thm}}
\newcommand{\elema}{\end{lema}}
\newcommand{\epro}{\end{pro}}
\newcommand{\ecoro}{\end{coro}}
\newcommand{\exmpl}{\end{xmpl}}
\newcommand{\ermrk}{\end{rmrk}}

\newcommand{\half}{\frac{1}{2}}

\newcommand{\B}{\mathcal B}

\newcommand{\EXP}{\textbf{e}}

\numberwithin{equation}{section}
\begin{document}
\begin{center}
{\bf{\large Characterization of unitary  processes with independent
and stationary
increments}}\\

\vspace{0.2in} {\large  Lingaraj Sahu {\footnote { Stat-Math Unit,
Indian Statistical Institute, Bangalore Centre, $8^{th}$ Mile,
Mysore Road,
 Bangalore-59, India.
 E-mail~:~lingaraj@isibang.ac.in }}  and
Kalyan B. Sinha {\footnote {  Jawaharlal Nehru Centre for Advanced
Scientific Research, Jakkur, Bangalore-64, India } } {\footnote {
Department of Mathematics, Indian Institute of Science,
Bangalore-12, India.\\ E-mail: kbs\_jaya@yahoo.co.in \\
 $^{1,2,3}$  The first and second authors acknowledge  support
 from National Board for Higher Mathematics, DAE  and
from  Bhatnagar Fellowship of  CSIR respectively.}}}
\end{center}

\begin{abstract}   This is  a continuation of the  earlier work
\cite{SSS} to  characterize  stationary unitary increment Gaussian
processes.  The earlier assumption of  uniform continuity is
replaced by weak continuity  and with a technical assumption on the
domain of the generator, unitary equivalence of the processes to the
solution of Hudson-Parthasarathy equation is proved.
\end{abstract}

\section{Introduction}
 In \cite{s1,s2}, by a co-algebraic
treatment,  Sch\"{u}rmann has proved that any weakly continuous
unitary stationary independent increment process on Hilbert space
$\mathbf h\otimes \mathcal H $ ( $ \mathbf h$  finite dimensional),
is unitarily equivalent to the solution of a Hudson-Parthasarathy
(HP) type quantum stochastic differential equation \cite{hp1}
 \be
\label{hpeqn00} dV_t=\sum_{\mu,\nu\ge 0}  V_t L_\nu^\mu
\Lambda_\mu^\nu(dt),~ V_0=1_{\mathbf h \otimes \Gamma} \ee where
$\Lambda_\mu^\nu$ are fundamental processes in the symmetric Fock
space  $ \Gamma(L^2(\mathbb R_+, \mathbf k))$ with respect to  a
fixed orthonormal basis (onb)  of  the noise  space $\mathbf k$  and
the coefficients $L_\nu^\mu~: \mu,~\nu~\ge 0$ are operators in the
initial Hilbert space  $\mathbf h$ given by
\begin{equation} \label{Lmunu}
L_\nu^\mu=\left\{ \begin{array} {lll}
 & G & \mbox{for}\
(\mu,\nu)=(0,0)\\
&L_j &  \mbox{for}\ (\mu,\nu)=(j,0)
\\
& -\sum_{j\ge 1}L_j^*W_k^j &
\mbox{for}\ (\mu,\nu)=(0,k)\\
& W_k^j- \delta_k^j 1_{\mathbf h} & \mbox{for}\ (\mu,\nu)=(j,k)
 \end{array}
 \right.
  \end{equation}
  ($\delta_k^j $ stands  for Dirac delta  function of $j$  and $k$) for some  operators   $G,~L_j$  in $\mathbf h$  and a unitary
operators  $W$  on $ \mathbf h \otimes  \mathbf k.$

For characterization of  Fock  adapted unitary evolution see
\cite{hl,ajl}  and references therein. In \cite{lw1,lw2}, by
extended semigroup methods,  Lindsay   and Wills have studied such
problems for Fock adapted contractive operator cocycles and
completely positive cocycles.

Recently in \cite{SSS} authors have studied the case of a unitary
stationary independent increment process on Hilbert space $\mathbf
h\otimes \mathcal H $  ( $ \mathbf h$  a separable Hilbert space),
with norm-continuous expectation semigroup  and showed its  unitary
equivalent to a Hudson-Parthasarathy    flow. Here  we are
interested in  unitary processes with  weakly continuous (not
necessarily  uniformly continuous )  expectation semigroup. Under
certain assumptions  on the domain of  the unbounded
 generators,  extending the ideas of \cite{SSS} we are
able to construct the noise space $\mathbf k$ and the operators
(unbounded) $G, L_j:\ge 1$   (see Proposition \ref{noise} and Lemma
\ref{H,Lj11}) such that the Hudson-Parthasarathy flow equation
(\ref{hpeqn00}) with coefficients (\ref{Lmunu})  (with $W$ being
identity operator), admits a unique unitary solution and the
solution is unitarily equivalent to the unitary process we started
with (see Theorem \ref{mainthm}).

\section{Notation and Preliminaries}

 We assume that  all the Hilbert spaces appearing in this article are complex
 separable with inner product  anti-linear
 in the first variable. For any Hilbert spaces $\mathcal H $ and $\mathcal K $ we denote the  Banach space  of bounded linear
operators from  $\mathcal H$ to
  $\mathcal K$ and trace class operators on $\mathcal H$ by $ \mathcal
 B(\mathcal H,\mathcal K)$  and
 $\mathcal B_1( \mathcal H)$
respectively. For a linear map (not necessarily bounded ) $T$  we
write its domain as  $\mathcal D(T).$ We shall denote the trace on
$\mathcal B_1( \mathcal H)$ by simply $Tr.$ The
 von Neumann algebra of bounded linear operators on  $\mathcal H$ is denoted by $B(\mathcal H).$
 The  Banach space   $\mathcal B_1( \mathcal H, \mathcal K)\equiv \{ \rho\in  \mathcal B( \mathcal H, \mathcal K)
 :  |\rho| :=\sqrt{\rho^*\rho }\in \mathcal B_1( \mathcal H) \}$  with norm (Ref. Page no. 47 in   \cite{gelshi})
\[\|\rho\|_1=\|  ~|\rho|~\|_{\mathcal B_1( \mathcal H)}= \sup\{ \sum_{k\ge 1} |\langle \phi_k, \rho \psi_k \rangle |
:\{\phi_k\},\{\psi_k\}  \}\] ( $\{\phi_k\},\{\psi_k\}$ varies over
orthonormal bases of $ \mathcal K $ and  $\mathcal H$ respectively )
is the predual of
 $\mathcal B( \mathcal K, \mathcal H).$  For an element $ x\in  \mathcal B( \mathcal K, \mathcal H),$ ~
 $  \mathcal B_1( \mathcal H, \mathcal K) \ni \rho \mapsto Tr(x\rho)$ defines an
  element  of the  dual
  Banach space  $\mathcal B_1( \mathcal H, \mathcal K)^*.$
 For a linear map $T$ on the Banach space $\mathcal B_1( \mathcal H, \mathcal K)$ the adjoint $T^*$ on the dual
  $\mathcal B( \mathcal K, \mathcal H)$ is given by
 $ Tr (T^* (x)\rho):=Tr(x  T(\rho)),~ \forall  x\in  \mathcal B( \mathcal K, \mathcal H),
 ~\rho \in \mathcal B_1( \mathcal H, \mathcal K). $\\

\noindent
 For any   $\xi \in \mathcal H \otimes
\mathcal K, h \in \mathcal H$  the map \[ \mathcal K \ni k
\mapsto\langle \xi, h\otimes k
 \rangle\] defines  a bounded  linear  functional on $\mathcal K$
  and thus  by  Riesz's  theorem
there exists  a unique vector  $\langle \langle h, \xi \rangle
\rangle$ in $ \mathcal K$  such that
 \be \label{partinn}
\langle~\langle \langle h, \xi \rangle \rangle,~k\rangle= \langle
\xi, h\otimes k\rangle, \forall k\in \mathcal K.\ee In other words
$\langle \langle h,\xi \rangle \rangle=F_{h}^*\xi$ where $F_{h}\in
\B(\mathcal K ,\mathcal H \otimes \mathcal K)$ is given by $F_{h}
k=h \otimes k.$

 Let $\mathbf h$ and $\mathcal H$ be two Hilbert spaces with some
  orthonormal bases $\{e_j:j \ge 1\}$ and
 $\{\zeta_j:j\ge 1\}$ respectively.
For $A\in \mathcal B( \mathbf h \otimes \mathcal H)$ and $u,v\in
\mathbf h$ we define a  linear operator $A(u,v)\in \mathcal B(
\mathcal H)$ by
 \[\langle \xi_1,A(u,v) \xi_2\rangle=\langle u\otimes \xi_1 ,A ~v\otimes \xi_2 \rangle,~\forall \xi_1,\xi_2 \in \mathcal
H\]  and read off the following properties (for a proof see Lemma
2.1 in  \cite{SSS}):
 \blema \label{Auv}
 Let $A,B \in \mathcal B( \mathbf h \otimes \mathcal H)$ then  for any $u,v,u_i$ and $v_i, i=1,2$ in $\mathbf h$
 \begin{description}
 \item  [(i)] $\|A(u,v)\| \le  \|A\|~\|u\|~ \|v\|$ and  $A(u,v)^*=A^*(v,u),$
 \item  [(ii)] $\mathbf h \times \mathbf h \mapsto A(\cdot~,\cdot)$ is  continuous bi-linear (anti-linear in first
 variable) mapping. If  $A(u,v)=B(u,v),~\forall u,v \in \mathbf h $ then $A=B,$
%\item  [(iii)] If a sequence of operator $\{A_n}_{n\ge 1}$ converges to $A$
%then $A_n(u,v)) $ converges to $A(u,v),~ \forall u,v \in \mathbf h$ in
%respective topology.
 \item  [(iii)]  $A(u_1,v_1)B(u_2,v_2)=[A(|v_1><u_2|\otimes 1_{\mathcal H})B](u_1,v_2),$
 \item  [(iv)] $AB(u,v)=\sum_{j \ge 1} A(u,e_j) B(e_j,v)$
 (strongly),
 \item  [(v)]  $0\le A(u,v)^* A(u,v) \le \|u\|^2 A^*A(v,v),$
 \item  [(vi)] $ \langle A(u,v)\xi_1, B(p,w)\xi_2 \rangle= \sum_{j\ge 1}\langle
 p \otimes \zeta_j, [B( |w><v|\otimes |\xi_2><\xi_1|)A^* u\otimes \zeta_j
\rangle\\
~~~~~~~~~~~~~~~~~~~~~~~~~~=\langle v\otimes \xi_1,~ [A^*(
|u><p|\otimes 1_{\mathcal H})B w \otimes \xi_2 \rangle .$
 \end{description}
 \elema

We also need to introduce  partial  trace  $Tr_{\mathcal H}$  which
is a linear map  from\\
 $\mathcal B_1( \mathbf h\otimes  \mathcal H)$  to $\mathcal B_1( \mathbf
 h)$ define  by, for $ B\in \mathcal B_1( \mathbf h\otimes  \mathcal H), $
\[ \langle u,Tr_{\mathcal H} (B) v\rangle:=\sum_{j\ge 1} \langle u\otimes \xi_j , B  v\otimes \xi_j
\rangle, \forall u,v\in \mathbf h.\] In particular,  for
$B=B_1\otimes B_2,  Tr_{\mathcal H} (B)= Tr (B_2)B_1.$

 For $A \in \mathcal B( \mathbf h \otimes \mathcal
H),\epsilon \in \mathbb Z_2=\{0,1\}$ we define operator
$A^{(\epsilon)}\in \mathcal B(\mathbf h\otimes \mathcal H)$ by
$A^{(\epsilon)}:= A$ if $\epsilon=0$  and $A^{(\epsilon)}:= A^*$ if
$\epsilon=1.$  For $1\le k\le n,$ we define a unitary exchange map
$P_{k,n}:\mathbf h^{\otimes n} \otimes \mathcal H\rightarrow \mathbf
h^{\otimes n} \otimes \mathcal H  $ by putting
 \[P_{k,n}(  u_1\otimes \cdots
  \otimes  u_n  \otimes  \xi ) :=u_1\otimes \cdots
  \otimes u_{k-1} \otimes u_{k+1}\cdots  \otimes  u_n \otimes  u_k \otimes \xi   \]
    on product vectors.
    %and then extending linearly.
%Since by definition range of $P_{k,n}$ is dense in  $\mathbf
%h^{\otimes n}\otimes \mathcal H,$ unitarity of $P_{k,n}$ follows
%from the fact that
%\[ \langle  P_{k,n}( \underbar{u}\otimes \xi), P_{k,n}( \underbar{v} \otimes \zeta ~ )
% \rangle =  \langle u_k \otimes  \xi, v_k
%  \otimes  \zeta \rangle \prod_{i\ne k}\langle  u_i,v_i  \rangle =\langle
 % \underbar{u}
 %  \otimes \xi,  \underbar{v} \otimes \zeta  \rangle.\]
Let $\underline{\epsilon}=(\epsilon_1,\epsilon_2,\cdots,
\epsilon_n)\in \mathbb Z_2^n.$  Consider the  ampliation of the
operator $A^{(\epsilon_k)}$ in $\mathcal B(\mathbf h^{\otimes n}
\otimes \mathcal H)$ given by
  \[ A^{(n,\epsilon_k)}:=P_{k,n}^* (1_{\mathbf h^{\otimes n-1}} \otimes A^{(\epsilon_k)})P_{k,n}.\]
   Now we define the operator   $A^{(\underline{\epsilon})}:=\prod
_{k=1}^n~A^{(n,\epsilon_k)}:=A^{(1,\epsilon_1)}\cdots
A^{(n,\epsilon_n)}$ in
  $\mathcal B(\mathbf h^{\otimes n} \otimes \mathcal H).$ Note
  that as
here, through out  this article, the product symbol  $\prod
_{k=1}^n$  stands for  product with the ordering   $1,2$ to $n.$
 For product vectors  $\underbar{u}, \underbar{v}\in
\mathbf h^{\otimes n}$  one can see  that
 \be \label{Akn} \prod_{i=1}^m A^{(n,\epsilon_i)} (\underbar{u},\underbar{v})
 =\prod_{i=1}^m A^{(\epsilon_i)}(u_i,v_i)   \prod_{i=m+1}^n \langle u_i,v_i\rangle
      \in  \mathcal B(\mathcal H).\ee
       When
$\underline{\epsilon}= \underline{0}\in \mathbb Z_2^n,$ for
simplicity we shall  write  $ A^{(n,k)} $ for
 $A^{(n,\epsilon_k)}$ and $A^{(n)}  $ for
$A^{(\underline{\epsilon})}.$

\subsection{Symmetric Fock Space  and Quantum Stochastic Calculus}

Let us briefly recall the
 fundamental integrator processes  of
quantum stochastic calculus  and the flow equation, introduced by
Hudson and Parthasarathy \cite{hp1}. For a Hilbert space $\mathbf k$
let us consider the symmetric Fock space $\Gamma=\Gamma(L^2(\mathbb
R_+,\mathbf k)).$  The exponential vector in  the Fock space,
associated with a vector $f\in L^2(\mathbb R_+, \mathbf k)$ is given
by
 \[\EXP(f)=\bigoplus_{n\ge 0}\frac{1}{\sqrt{n!}}f^{(n)},\]
 where $ f^{(n)}=\underbrace{f \otimes f  \otimes
 \cdots \otimes  f}_{n-copies} $ for $ n>0 $
and by convention  $f^{(0)}=1 .$ The exponential vector
$\textbf{e}(0)  $ is called    the vacuum vector. For any subset
 $M$ of $L^2(\mathbb R_+, \mathbf k)$  we shall  write  $\mathcal E(M)$  for  the  subspace  spanned  by
 $\{\EXP(f):f\in M\}.$ For an interval $\Delta$
 of $\mathbb R_+,$ let $\Gamma_\Delta$ be the symmetric Fock space over
  the Hilbert space $L^2(\Delta, \mathbf k) \cong$ the range of the multiplication operator
   $1_\Delta$
  on $L^2(\mathbb R_+, \mathbf k).$
For  $0\le s\le t<\infty,$
 the Hilbert space $\Gamma$ decompose as
 $\Gamma_{s]} \otimes  \Gamma_{(s,t]} \otimes
\Gamma_{[t}$ respectively, here we have abbreviated $[0,s]$ by $s]$
and $(t, \infty)$ by $[t,$  and for any $f\in L^2(\mathbb R_+,
\mathbf k)$  the exponential vector $\EXP(f)= \EXP(f_{s]}) \otimes
\EXP( f_{(s,t]}) \otimes \EXP(f_{[t})$ where $f_\Delta= 1_\Delta f.$

 \noindent  Let us consider the
Hudson-Parthasarathy  (HP) flow equation on $\mathbf h \otimes
\Gamma(L^2(\mathbb R_+, \mathbf k))$:

\be \label{hpeqn} V_{s,t}=1_{\mathbf h \otimes
\Gamma}+\sum_{\mu,\nu\ge 0}  \int_s^t V_{s,\tau} L_\nu^\mu
\Lambda_\mu^\nu(d\tau).\ee Here the coefficients $L_\nu^\mu~:
\mu,~\nu~\ge 0$ are operators in $\mathbf h$  (not necessarily
bounded) and $\Lambda_\mu^\nu$ are fundamental processes with
respect to a fixed orthonormal basis $\{E_j: j\ge 1\}$ of $\mathbf
k:$

 \begin{equation}
 \Lambda_\nu^\mu (t)=\left\{
 \begin{array} {lll} & t ~1_{\mathbf h \otimes
\Gamma} & \mbox{for}
(\mu,\nu)=(0,0) \\
& a(1_{[0,t]}\otimes E_j)&  \mbox{for} (\mu,\nu)=(j,0)
\\
& a^\dag(1_{[0,t]}\otimes E_k)& \mbox{for}(\mu,\nu)=(0,k)
\\
& \Lambda (1_{[0,t]}\otimes |E_k><E_j|)& \mbox{for} (\mu,\nu)=(j,k).
 \end{array}
 \right.
  \end{equation}
The fundamental processes  $a, a^\dag$  and $\Lambda$  are called
annihilation, creation  and  conservation   respectively (for their
definition  and   detail  about quantum stochastic calculus see
\cite{krp, gs}).

\section {Unitary processes  with stationary and independent  increments}

Let $\{U_{s,t}: 0\le s\le t<\infty\}$ be a family of unitary
operators in $ \mathcal B( \mathbf h \otimes \mathcal H)$ and
$\Omega $ be a fixed unit vector in $\mathcal H.$ We shall write
$U_t:=U_{0,t}$ for simplicity. Let us consider the family of unitary
operators  $\{ U_{s,t}^{(\epsilon)}\}$ in $ \mathcal B( \mathbf
h\otimes \mathcal H)$  for $\epsilon \in \mathbb Z_2$ given by $
U_{s,t}^{(\epsilon)}= U_{s,t}$ if $\epsilon=0, U_{s,t}^{(\epsilon)}=
U_{s,t}^*$ if $\epsilon=1.$ As in previous section,  for $n\ge 1 ,
\underline{\epsilon}\in \mathbb Z_2^n$ fixed  and  $1\le k\le n, $
we define the families  of operators $\{U_{s,t}^{(n,\epsilon_k)}\}$
and $\{U_{s,t}^{(\underline{\epsilon})}\}$ in $ \mathcal B( \mathbf
h^{\otimes n} \otimes \mathcal H).$    By identity  (\ref{Akn}) we
have, for product vectors  $\underbar{u},\underbar{v}\in \mathbf
h^{\otimes n}$ and $\underline{\epsilon}\in \mathbb Z_2^n,$
\[
U_{s,t}^{(\underline{\epsilon})}(\underbar{u},\underbar{v})
=\prod_{i=1}^n U_{s,t}^{(\epsilon_i)}(u_i,v_i).\]

\noindent  Furthermore,  for $ \underbar{s}=(s_1,s_2, \cdots, s_n),
\underbar{t}=(t_1,t_2, \cdots, t_n)$ $:~ 0 \le s_1\le t_1\le s_2\le
\ldots \le  s_n\le t_n< \infty,$ we define $
U_{\underbar{s},\underbar{t}}^{(\underline{\epsilon})}\in \mathcal
B(\mathbf h^{\otimes n}\otimes \mathcal H)$ by setting \be
\label{U-underbar-st}
U_{\underbar{s},\underbar{t}}^{(\underline{\epsilon})}:=\prod_{k=1}^n
U_{s_k,t_k}^{(n,\epsilon_k)}.\ee
 Then for $\underbar{u}=\otimes_{k=1}^n u_k,
\underbar{v}=\otimes_{k=1}^n v_k\in \mathbf h^{\otimes n} $ we have
\[U_{\underbar{s},\underbar{t}}^{(\underline{\epsilon})}(\underbar{u},
\underbar{v}) =\prod_{k=1}^n U_{s_k,t_k}^{(\epsilon_k)} (u_k,v_k)
.\] When $\underline{\epsilon}=\underline{0},$ we  write
$U_{\underbar{s},\underbar{t}}$ for
$U_{\underbar{s},\underbar{t}}^{(\underline{\epsilon})}.$ For
$\alpha, \beta \ge 0, \underbar{s}=(s_1,s_2, \cdots, s_n),
\underbar{t}=(t_1,t_2, \cdots, t_n)$ we write $\alpha\le
\underbar{s},\underbar{t}\le \beta$  if $ \alpha\le s_1\le t_1\le
s_2 \le \ldots \le s_n\le  t_n\le \beta.$

\noindent  We assume the following on the family of unitary
$\{U_{s,t}\in \B(\mathbf h \otimes \mathcal H)\}.$
  \\
 {\bf Assumption A }
\begin{description}
\item [A1] {\bf (Evolution)} For any $   0\le r\le s\le t<\infty, ~ U_{r,s}U_{s,t}=U_{r,t}.$
\item [A2] {\bf (Independence of increments)} For any ~ $0\le s_i\le t_i<\infty ~:
 ~ i=1,2$ such that $[s_1,t_1)\cap  [s_2,t_2)=\emptyset$\\
 (i)~$
U_{s_1,t_1}(u_1,v_1) $ commutes with  $U_{s_2,t_2}(u_2,v_2) $
 and  $U_{s_2,t_2}^*(u_2,v_2) $  for every  $u_i,v_i \in \mathbf
 h.$\\
(ii)~  For  $s_1\le \underbar{a},\underbar{b}\le  t_1,~~s_2\le
\underbar{q},\underbar{r}\le  t_2$ and $\underbar{u},\underbar{v}\in
\mathbf h^{\otimes n},~
 \underbar{p},\underbar{w}\in
\mathbf h^{\otimes m},\underline{\epsilon}\in \mathbb
Z_2^n,\underline{\epsilon}^\prime \in \mathbb Z_2^m $
\[\langle \Omega,  U_{\underbar{a},\underbar{b}}^{(\underline{\epsilon})}(\underbar{u},
\underbar{v})
U_{\underbar{q},\underbar{r}}^{(\underline{\epsilon}^\prime)}(\underbar{p},
\underbar{w}) \Omega \rangle = \langle \Omega,
U_{\underbar{a},\underbar{b}}^{(\underline{\epsilon})}(\underbar{u},
\underbar{v})\Omega \rangle
 \langle \Omega,U_{\underbar{q},\underbar{r}}^{(\underline{\epsilon}^\prime)}(\underbar{p},
\underbar{w})  \Omega \rangle.\]
\item [A3]  {\bf (Stationarity of increments)} For any $0\le s\le t<\infty$
 and  $ \underbar{u},\underbar{v}\in \mathbf h^{\otimes n} ,
 \underline{\epsilon}\in \mathbb Z_2^n$
 \[\langle \Omega,  U_{s,t}^{(\underline{\epsilon})}(\underbar{u},\underbar{v}
 ) \Omega \rangle
=  \langle \Omega,
U_{t-s}^{(\underline{\epsilon})}(\underbar{u},\underbar{v}) \Omega
\rangle.\]
\end{description}
\noindent {\bf Assumption  $B^\prime$ ~(Weak / Strong continuity)}
\[\lim_{t \rightarrow 0 }~ \langle \Omega,
 (U_t-1)(u,v) \Omega \rangle =0,~\forall u,v \in \mathbf h.\]

\brmrk The  assumption  {\bf  $B^\prime$} is an weakening of the
{\bf assumption  B} in \cite{SSS}. \ermrk

\noindent  As in \cite{SSS} we also assume  the following
simplifying conditions.
\begin{description}
\item [Assumption  C   (Gaussian condition) ]
 For any $u_i,v_i\in \mathbf h,$ \\ $\epsilon_i\in \mathbb
Z_2: i=1,2,3$ \be \lim_{t\rightarrow 0}\frac{1}{t} \langle
\Omega,~(U_t^{(\epsilon_1)}-1)(u_1,v_1)(U_t^{(\epsilon_2)}-1)(u_2,v_2)
(U_t^{(\epsilon_3)}-1)(u_3,v_3) ~\Omega \rangle=0. \ee
\end{description}

\begin{description}
\item[Assumption D ]  {\bf (Minimality)}
 The set
$ \mathcal S_0=\{
U_{\underbar{s},\underbar{t}}(\underbar{u},\underbar{v})
\Omega:=U_{s_1,t_1}(u_1,v_1)\cdots U_{s_n,t_n}(u_n,v_n)\Omega :
\underbar{s}=(s_1,s_2, \cdots, s_n),$ $ \underbar{t}=(t_1,t_2,
\cdots, t_n)$ $:~ 0 \le s_1\le t_1\le  s_2 \cdots, s_n\le  t_n<
\infty,n\ge 1, \underbar{u}=\otimes_{i=1}^n
u_i,\underbar{v}=\otimes_{i=1}^n v_i~\mbox{with}~ u_i,v_i\in \mathbf
h\}$ is total in $\mathcal H.$
\end{description}

\brmrk  The  assumption  {\bf D  }  is  not really a restriction,
one
 can as well work with replacing $\mathcal  H $  by span closure  of $\mathcal
 S_0.$
 \ermrk
 \brmrk \label{DS_0}  For  any  dense set  $\mathcal D\subseteq\mathbf h,  \mathcal S_0$ will be  still  total
 if we  restrict $ u_i,v_i\in
 \mathcal D$  in   the  assumption {\bf  D}.
\ermrk

\subsection{Expectation Semigroups} Let us look at the various
semigroups associated with the evolution $\{U_{s,t}\}.$

\noindent For any fixed  $n\ge 1,$ we define a family of  operators
$\{T_t^{(n)}\}$ on $\mathbf h ^{\otimes n}$ by setting
\[\langle \phi,T_t^{(n)} ~\psi \rangle:= \langle \Omega,
  U_t^{(n)}(\phi,\psi)~ \Omega \rangle,~ \forall
  \phi,\psi \in { \mathbf h}^{\otimes n}.\]
Then in particular  for product vectors $\underbar{u}=
\otimes_{i=1}^n u_i,~ \underbar{v}=\otimes_{i=1}^n v_i\in\mathbf h
^{\otimes n}$
\[\langle \underbar{u},T_t^{(n)} ~\underbar{v}   \rangle= \langle \Omega,
U_t^{(n)}(\underbar{u},\underbar{v})~ \Omega \rangle = \langle
\Omega,   U_t(u_1,v_1) U_t(u_2,v_2)\cdots U_t(u_n,v_n)~ \Omega
\rangle.\] We shall write $T_t$  for  $T_t^{(1)}.$

\bpro \label{Ttsemi}
  Under the  assumption {\bf  A}  and {\bf $B^\prime$} the
 $\{T_t^{(n)} \}$ for each $n\ge 1$ is a   strongly continuous contractive semigroup on $\mathbf h^{\otimes n}$.
\epro

We need a Lemma for the proof of this proposition. That $T_t^{(n)}$
is a semigroup  follows exactly as in the proof of  Lemma 6.1  in
\cite{SSS}  which as well as that of  following Lemma we omit.

 \blema \label{UTknc0}
\begin{description}
\item[(i)] For $1\le k\le n,$
 \be
\langle \Omega ,U_t^{(n,k)} (\underbar{p}, \underbar{w}) \Omega
\rangle =\langle \underbar{p}, 1_{{\mathbf h} ^{(\otimes k-1)}}
\otimes T_t \otimes 1_{{\mathbf h}^{(\otimes n-k)}} \underbar{w}
\rangle, ~\forall \underbar{p},\underbar{w} \in { \mathbf
h}^{\otimes n}. \ee We shall denote this ampliation $ 1_{{\mathbf
h}^{(\otimes k-1)}} \otimes T_t \otimes 1_{{\mathbf h}^{(\otimes
n-k)}}$  by $T_t^{(n,k)}  $.
\item[(ii)] For  any $1\le m\le n,~ \underbar{p},\underbar{w}
\in { \mathbf h}^{\otimes n},$
\[
 \langle \Omega ,(\prod_{k=1}^m U_t^{(n,k)})  (\underbar{p}, \underbar{w}) \Omega \rangle
=\langle  \underbar{p}, T_t^{(m)}\otimes 1_{{\mathbf h}^{(\otimes
n-m)} }~ \underbar{w} \rangle. \]
\item[(iii)]
 For any $\phi\in { \mathbf h}^{\otimes n},$
\bean
&&\|(U_t^{(n,k)}-1) \phi\otimes \Omega \|^2\\
&& = \langle (1-T_t^{(n,k)}) \phi,\phi \rangle
+\langle \phi,(1-T_t^{(n,k)}) \phi\rangle \\
&&\le 2\|(1-T_t)\phi\|~\|\phi\|. \eean
\item[(iv)] For any $\phi\in { \mathbf h}^{\otimes n},$
 \bean
&& \|(U_t^{(n)} -1) \phi\otimes \Omega \|^2\\
&& = \langle (1-T_t^{(n)}) \phi,\phi\rangle
+\langle \phi,(1-T_t^{(n)}) \phi\rangle \\
&&\le 2\|(1-T_t^{(n)}) \phi \|~\|\phi\|.\eean
\item[(v)] For any $v\in \mathbf h$
\be \sum_{m\ge 1} \|(U_t-1)(e_m,v)\Omega\|^2=2 Re \langle v ,
(1-T_t) v \rangle \le 2 \|v\|~  \|(T_t-1)v\|.\ee

\end{description}
\elema
 \noindent  {\bf Proof of the Proposition \ref{Ttsemi}  :}\\
 \noindent The  assumption {\bf
$B^\prime$} and definition of $T_t$ implies that the semigroup of
contractions $ \{T_t\}$ on $\mathbf h$ is weakly and hence strongly
continuous. To apply induction let us assume that for some $m\ge 1,$
the contractive semigroups $\{T_t^{(n)}\} $ are strongly continuous
for all $1\le n \le m-1.$ Now let us consider the following, for any
$ \phi,\psi\in {\mathbf h}^{\otimes m},$ \bean
&& \langle \phi \otimes \Omega,   (U_t^{(m)} -1) \psi \otimes \Omega \rangle\\
&& =\langle \phi\otimes \Omega,   \left (
[\prod_{k=1}^{m-1}U_t^{(m,k)}]
[U_t^{(m,m)}]-1 \right )\psi\otimes  \Omega \rangle\\
&& =\langle  [\prod_{k=1}^{m-1}U_t^{(m,k)}] ^*\phi\otimes \Omega,
\left ([U_t^{(m,m)}]-1 \right )\psi\otimes  \Omega \rangle\\
&& ~~~+\langle \phi\otimes \Omega,   \left (
[\prod_{k=1}^{m-1}U_t^{(m,k)}]-1 \right )\psi\otimes  \Omega
\rangle.\eean Taking absolute value, by Lemma \ref{UTknc0}~we get
\bean
&& |\langle \phi,~ (T_t^{(m)}-1_{{\mathbf h}^{\otimes m }})\psi \rangle|\\
&& \le \| \phi\|~     \sqrt{2~\|\psi\|~\|[(1_{{\mathbf h}^{\otimes
m-1 }} \otimes T_t)- 1_{{\mathbf h}^{\otimes m }}] \psi \|}
 + |\langle \phi,   \left ( [T_t^{(m-1)} \otimes 1_{\mathbf h}]-
 1_{{\mathbf h}^{\otimes m }} \right )\psi \rangle|\\
&&\le \|\phi\|   \sqrt{2~\|\psi\|~\|  [1_{{\mathbf h}^{\otimes m-1
}} \otimes ( T_t- 1_{\mathbf h})] \psi \|} + \|\phi\| ~ \|(
[T_t^{(m-1)}-1_{{\mathbf h}^{\otimes m-1 }}]  \otimes 1_{\mathbf h}
) \psi\| . \eean So strong continuity of $T_t^{(m-1)}$ and $T_t$
implies $T_t^{(m)}$ is strongly continuous. \qed

\noindent Let us denote the  generator of the semigroup $T_t^{(n)}$
 by $G^{(n)}$ and for $n=1$ by  $G$ with   domain $\mathcal D(G).$
 % Since $G$ is the
 % generator  of a contraction semigroup, $\mathcal D(G^2)$  is a
%core for $G.$

 \blema \label{4Ut-111} Under the  assumption {\bf  C} we have the
 following.\\
{\bf (i)} For any $n\ge 3,~\underbar{u},\underbar{v} \in \mathbf
h^{\otimes n},\underline{\epsilon}\in \mathbb Z_2^n$
 \be  \lim_{t\rightarrow
0}\frac{1}{t} \langle \Omega,~ (U_t^{(\epsilon_1)}-1)(u_1,v_1)\cdots
(U_t^{(\epsilon_n)}-1)(u_n,v_n)~\Omega \rangle=0 .\ee \noindent {\bf
(ii)} For vectors  $u\in \mathbf h, v\in \mathcal D(G),$
   product vectors $ \underbar{p},\underbar{w}\in \mathbf h^{\otimes n}$ and $\epsilon
\in \mathbb Z_2, \underline{\epsilon^\prime} \in \mathbb Z_2^n$
 \be \label{UtUt*inner}
\lim_{t\rightarrow 0}\frac{1}{t} \langle (U_t-1)^{(\epsilon)}(u,v)
~\Omega,
(U_t^{(\underline{\epsilon^\prime})}-1)(\underbar{p},\underbar{w})~\Omega
\rangle\ee
\[= (-1)^\epsilon
\lim_{t\rightarrow 0}\frac{1}{t} \langle (U_t-1)(u,v) ~\Omega,
(U_t^{(\underline{\epsilon^\prime})}-1)(\underbar{p},\underbar{w})~\Omega
\rangle.\] \elema
\begin{proof} {\bf (i)} The  proof is identical  to that of Lemma 6.7
in
\cite{SSS}.\\
\noindent {\bf (ii)}  For $\epsilon =0$ nothing to prove.   To see
this for $\epsilon =1$ consider the following \bea \label{**}
&&\lim_{t\rightarrow 0}\frac{1}{t}\langle (U_t+U_t^*-2)(u,v) \Omega,
(U_t^{(\underline{\epsilon^\prime})}-1)(\underbar{p},\underbar{w})~\Omega\rangle \\
&&=-\lim_{t\rightarrow 0}\frac{1}{t}\langle[ (U_t^*-1)(U_t-1)](u,v)
\Omega,
(U_t^{(\underline{\epsilon^\prime})}-1)(\underbar{p},\underbar{w})~\Omega\rangle
{\nonumber}
\\
&&=-\lim_{t\rightarrow 0}\frac{1}{t}
 \sum_{m\ge 1} \langle   (U_t-1)(e_m,v) \Omega,
  (U_t-1)(e_m,u)(U_t^{(\underline{\epsilon^\prime})}-1)
  (\underbar{p},\underbar{w})~\Omega\rangle. {\nonumber}
\eea

\noindent That this limit  vanishes   can be seen from the following
 \bean &&|
\frac{1}{t}
 \sum_{m\ge 1} \langle   (U_t-1)(e_m,v) \Omega,
  (U_t-1)(e_m,u)(U_t^{(\underline{\epsilon^\prime})}-1)
  (\underbar{p},\underbar{w})~\Omega\rangle|^2\\
  &&\le
 \sum_{m\ge 1} \frac{1}{t} \|(U_t-1)(e_m,v) \Omega\|^2
~~ \sum_{m\ge 1} \frac{1}{t}
\|(U_t-1)(e_m,u)(U_t^{(\underline{\epsilon^\prime})}-1)
  (\underbar{p},\underbar{w})~\Omega\|^2.\eean
By Lemma  \ref{UTknc0} (v) and Lemma  \ref{Auv} (iv) the above
quantity is equal to \bean && 2 Re \langle v, \frac{1-T_t}{t}
v\rangle  \frac{1}{t} \langle
(U_t^{(\underline{\epsilon^\prime})}-1)
  (\underbar{p},\underbar{w})~\Omega,
 [(U_t^*-1)(U_t-1)](u,u) (U_t^{(\underline{\epsilon^\prime})}-1)
  (\underbar{p},\underbar{w})~\Omega  \rangle\\
&&\le 2 Re \langle v, \frac{1-T_t}{t} v\rangle  \frac{1}{t} \langle
(U_t^{(\underline{\epsilon^\prime})}-1)
  (\underbar{p},\underbar{w})~\Omega,
 (2-U_t^*-U_t)(u,u) (U_t^{(\underline{\epsilon^\prime})}-1)
  (\underbar{p},\underbar{w})~\Omega  \rangle
\eean Therefore, since $Re \langle v, \frac{1-T_t}{t} v\rangle $ is
uniformly  bounded in $t$  as  $T_t$ is strongly   continuous  and
$v\in \mathcal D(G)$, by  assumption {\bf C} we get
\[\lim_{t\rightarrow 0}\frac{1}{t}
 \sum_{m\ge 1} \langle   (U_t-1)(e_m,u) \Omega,
  (U_t-1)(e_m,v)(U_t^{(\underline{\epsilon^\prime})}-1)
  (\underbar{p},\underbar{w})~\Omega\rangle=0.\]
  Thus (\ref{UtUt*inner}) follows.
  \end{proof}

\noindent For  vectors  $u,p\in \mathbf h$ and $v,w\in \mathcal
D(G),$ the identity (\ref{UtUt*inner}) gives \be \label{UtUt*10}
\lim_{t\rightarrow 0}\frac{1}{t} \langle (U_t-1)^{(\epsilon)}(u,v)
~\Omega, (U_t-1)^{\epsilon^\prime}(p,w) ~\Omega\rangle\ee
\[= (-1)^{\epsilon+\epsilon^\prime}
\lim_{t\rightarrow 0}\frac{1}{t} \langle (U_t-1)(u,v) ~\Omega,
(U_t-1)(p,w) ~\Omega\rangle.\]

\noindent For $m,n\ge 1,$ we define a family of operators $\{Z_t^{(
m,n)}:t\ge 0\}$ on the Banach space $\mathcal B_1( \mathbf
h^{\otimes m}, \mathbf h^{\otimes n})$ by

\[Z_t^{(   m,n)} \rho=  Tr_{\mathcal H}[ U_t^{(n)} ( \rho \otimes |\Omega ><  \Omega|) \
(U_t^{(m)})^*],~\rho \in \mathcal B_1( \mathbf h^{\otimes m}, \mathbf h^{\otimes n}) .\]
Then in particular for product vectors $\underbar{u},\underbar{v}
\in \mathbf h^{\otimes m},\underbar{p},\underbar{w} \in \mathbf
h^{\otimes n}.$ \be \label{Ztsimple} \langle \underbar{p} ,Z_t^{(
m,n)} (|\underbar{w}><\underbar{v}|) \underbar{u}  \rangle:= \langle
 U_t^{(m)}(\underbar{u},\underbar{v}) \Omega,~ U_t^{(n)}(\underbar{p},\underbar{w})~\Omega
\rangle.\ee

\blema The above family  $\{Z_t^{( m,n)}\}$ is a  semigroup of
contractive maps  on $\mathcal B_1( \mathbf h^{\otimes m}, \mathbf
h^{\otimes n}).$  Furthermore  assumption {\bf $B^\prime$ } implies
$\{Z_t^{( m,n)}\}$ is strongly continuous in the $\mathcal B_1$
topology . \elema
\begin{proof}
For $\rho\in \mathcal B_1( \mathbf h^{\otimes m}, \mathbf h^{\otimes
n})$ \bean
&&\|Z_t^{(m,n)} \rho\|_1=  \|Tr_{\mathcal H}[ U_t^{(n)} ( \rho \otimes |\Omega ><  \Omega|) (U_t^{(m)})^*]\|_1\\
&& = \sup_{\phi^{(l)}~ onb~ of~ \mathbf h ^{\otimes l}~:~ l=m,n }
~\sum_{k\ge 1} |\langle \phi_k^{(n)} ,Tr_{\mathcal H}[ U_t^{(n)}(
\rho \otimes
|\Omega ><  \Omega|) (U_t^{(m)})^*] \phi_k^{(m)} \rangle|\\
&& \le  \sup_{\phi^{(l)}} \sum_{j,k\ge 1} |\langle
\phi_k^{(n)}\otimes \zeta_j ,
 U_t^{(n)} ( \rho \otimes |\Omega ><  \Omega|) (U_t^{(m)})^* \phi_k^{(m)}
 \otimes \zeta_j \rangle|\\
&& \le \|U_t^{(n)} ( \rho \otimes |\Omega >< \Omega|)
(U_t^{(m)})^*\|_1. \eean Since for any $l\ge 1,~\{U_t^{(l)}\}$ is a
family of unitary operators
\[\|Z_t^{(m,n)} \rho\|_1 =  \|\rho \otimes |\Omega ><  \Omega|\|_1 =  \|\rho\|_1.\]
%and thus contractivity of $Z_t^{(m,n)}.$
Proof of semigroup property of  $\{Z_t^{(m,n)}\}$  is  same as in
Lemma 6.4  \cite{SSS}. \noindent In order to prove strong continuity
$ Z_t^{( m,n)},$ it is suffices to prove the same  for rank one
operator $\rho=|\underbar{w}><\underbar{v}|,
~\underbar{v},\underbar{w}$ product vectors in $ \mathbf h^{\otimes
m}$  and $ \mathbf h^{\otimes n}$  respectively. We have \bean
&&\|(Z_t^{(   m,n)}-1) (|\underbar{w}><\underbar{v}|)\|_1\\
&& = \sup_{\phi^{(l)}~ onb~ of~ \mathbf h ^{\otimes l}~:~
l=m,n}\sum_{k\ge 1} |\langle \phi_k^{(n)} ,
(Z_t^{(   m,n)} -1) (|\underbar{w}><\underbar{v}|)  \phi_k^{(m)} \rangle|\\
&& = \sup_{\phi^{(l)}}
 \sum_{k\ge 1} |\langle U_t^{(m)} ( \phi_k^{(m)},\underbar{v}) \Omega    ,
 U_t^{(n)}  ( \phi_k^{(n)}, \underbar{w} ) \Omega  \rangle-\overline{\langle
  \phi_k^{(m)},\underbar{v} \rangle} \langle \phi_k^{(n)}, \underbar{w}   \rangle|\\
&& \le \sup_{\phi^{(l)}} \sum_{k\ge 1} |\langle (U_t^{(m)}-1)
 ( \phi_k^{(m)},\underbar{v}) \Omega    , U_t^{(n)}  ( \phi_k^{(n)}, \underbar{w} )
 \Omega  \rangle|\\
&& +\sup_{\phi^{(l)}}  \sum_{k\ge 1} |\overline{\langle
\phi_k^{(m)},
\underbar{v} \rangle}  \langle  \Omega    , (U_t^{(n)}-1)  ( \phi_k^{(n)}, \underbar{w} ) \Omega  \rangle|\\
&& \le \sup_{\phi^{(l)}} \left [\sum_{k\ge 1} \| (U_t^{(m)}-1) (
\phi_k^{(m)},\underbar{v}) \Omega \|^2 \right ]^{\frac{1}{2}}
 \left [\sum_{k\ge 1} \|U_t^{(n)}  ( \phi_k^{(n)}, \underbar{w} ) \Omega \|^2 \right ]^{\frac{1}{2}} \\
&& +\sup_{\phi^{(l)}} \left [\sum_{k\ge 1} |\langle
\phi_k^{(m)},\underbar{v} \rangle|^2 \right ]^{\frac{1}{2}} \left
[\sum_{k\ge 1} \|(U_t^{(n)}-1)  ( \phi_k^{(n)}, \underbar{w} )
\Omega \|^2 \right ]^{\frac{1}{2}}. \eean Hence by Lemma
\ref{UTknc0} \bean
&&\|(Z_t^{(   m,n)}-1) (|\underbar{w}><\underbar{v}|)\|_1\\
&&\le \|\underbar{w}\| \sqrt{2 ~\|(T_t^{(m)}-1)\underbar{v}\|}+
~\|\underbar{v}\|\sqrt{2 \|(T_t^{(n)}-1)\underbar{w}\|}. \eean
\noindent  Thus by  strong   continuity of the semigroup $
T_t^{(m)}$ and $T_t^{(n)},$  and the density of the finite rank
vectors in $\mathcal B_1( \mathbf h^{\otimes m}, \mathbf h^{\otimes
n}) $ the contractive semigroup ${Z_t^{( m,n)}}$ is a strongly
continuous on $\mathcal B_1( \mathbf h^{\otimes m}, \mathbf
h^{\otimes n}).$
\end{proof}

 We shall  denote
the   generator of  the semigroup $Z_t^{( m,n)}$ by $\mathcal L^{(
m,n)}.$ For $n\ge 1$ we shall write $Z_t^{(n)}$ for   the semigroup
$Z_t^{( n,n)}$ on the Banach space $\mathcal B_1(\mathbf h^{\otimes
n})$ with denoting its generator by $\mathcal L^{(n)}$ for
simplicity. Moreover,  we denote the semigroup $Z_t^{(1)}$ and its
generator  $\mathcal L^{(1)}$ by just $Z_t$  and $\mathcal L$
respectively.

 \blema  \label{Zn+11} For any $n\ge 1 ,~Z_t^{(n)}$ is a
positive  trace preserving semigroup.
 \elema
\begin{proof}
Positivity  follows from the following, for any
$\underbar{u},\underbar{v} \in \mathbf h ^{\otimes n}$ \bean
&&\langle \underbar{u} ,Z_t^{(n)} (|\underbar{v}><\underbar{v}|)
 \underbar{u}  \rangle = \|  U_t^{(n)}(\underbar{u},\underbar{v})
 \Omega\|^2.
\eean By definition we have
%and  Lemma \ref{Auv}
 \bean
&& Tr[ Z_t^{(n)} (|\underbar{u}><\underbar{v}|)]= \sum_{k}\langle \underbar{e}_k ,Z_t^{(n)} (|\underbar{u}><\underbar{v}|) \underbar{e}_k  \rangle\\
&&=\sum_{k}\langle U_t^{(n)}(\underbar{e}_k,\underbar{v})\Omega ,U_t^{(n)}(\underbar{e}_k,\underbar{u}) \Omega \rangle\\
&&=\langle \Omega ,(U_t^{(n)})^*
U_t^{(n)}(\underbar{v},\underbar{u}) \Omega \rangle. \eean Since
$U_t^{(n)}$ is  unitary, we get \be  \label{tr-zero11} Tr[ Z_t^{(n)}
(|\underbar{u}><\underbar{v}|)]= \langle  \underbar{v}, \underbar{u}
\rangle= Tr (|\underbar{u}><\underbar{v}|).\ee
\end{proof}

\noindent Let us define a family  $\{Y_t:t\ge 0\}$  of positive
contractions on $\B_1(\mathbf h)$ by $Y_t(\rho):=T_t~\rho~ T_t^*,~~
\forall \rho \in \mathcal B_1(\mathbf h).$ Since $T_t$  is a $C_0$-
semigroup  of  contraction operators  on $\B(\mathbf h)$ it can be
seen that $Y_t$ is a contractive  $C_0$-semigroup  on $\B_1(\mathbf
h).$ It can also be seen that \cite{gs} the generator
$\widetilde{\mathcal L}$ of $Y_t$ satisfy
\[\widetilde{\mathcal L}(\rho)=G^* \rho+ \rho G,~~\forall \rho
\in \mathcal D_0\equiv \{(1-G)^{-1}\sigma  (1-G^*)^{-1} : \sigma \in
\mathcal B_1(\mathbf h)\}\] and    $\mathcal D_0$ is a core for
$\widetilde{\mathcal L}.$  If we define   the subspace $\mathcal
N_0\equiv Span \{|u><v|, u,v \in \mathcal D(G)\}$  of $\mathcal
B_1(\mathbf h),$  then it is clear that $ \mathcal N_0$ is dense in
$\mathcal B_1(\mathbf h)$  and contained in  $\mathcal D_0.$

\noindent  We also need another class of semigroup. For $m,n\ge 1$
we define a family of maps $F^{(m,n)}_t$ on the Banach space
$\mathcal B_1( \mathbf h^{\otimes m}, \mathbf h^{\otimes n})$ by
 \be F^{(m,n)}_t \rho=  Tr_{\mathcal H}[
(U_t^{(n)})^*
 ( \rho \otimes |\Omega ><  \Omega|) U_t^{(m)}], ~\forall  \rho
 \in  \mathcal B_1( \mathbf h^{\otimes m}, \mathbf h^{\otimes n})
 \ee

\noindent So in particular for product vectors
$\underbar{u},\underbar{v} \in \mathbf h^{\otimes m}$ and
$\underbar{p},\underbar{w} \in \mathbf
h^{\otimes n},$ we have that \\
  $\langle
\underbar{p} ,F^{(m,n)}_t (|\underbar{w}><\underbar{v}|)
\underbar{u}
  \rangle
=\langle   ( U_t^{(m)})^*(\underbar{u},\underbar{v}) \Omega,~ (
U_t^{(n)})^* (\underbar{p},\underbar{w})~\Omega \rangle. $ \blema
For any $m,n\ge 1, ~\{F^{(m,n)}_t:t\ge 0\} $ is a strongly
continuous contractive semigroup on $\mathcal B_1( \mathbf
h^{\otimes m}, \mathbf h^{\otimes n})$. \elema
\begin{proof} The proof is same as for the semigroup  ${Z_t^{( m,n)}}.$
\end{proof}
For $n= 1,$ we shall write $F_t$ for   the semigroup $F_t^{( 1,1)}$
on the Banach space $\mathcal B_1(\mathbf h)$ and  shall denote its
generator  by $\mathcal L^\prime. $

\section{Construction of noise space}
\noindent Let
$M_0:=\{(\underbar{u},\underbar{v},\underline{\epsilon}):
~\underbar{u}=\otimes_{i=1}^n u_i, \underbar{v}=\otimes_{i=1}^n v_i,
u_i\in \mathbf h, v_i \in   \mathcal D(G),
\underline{\epsilon}=(\epsilon_1,\cdots, \epsilon_n)\in \mathbb
Z_2^n, ~ n\ge 1\}$  and consider the relation $``\sim"$ on $M_0$ as
defined in \cite{SSS} :
$(\underbar{u},\underbar{v},\underline{\epsilon})\sim
(\underbar{p},\underbar{w},\underline{\epsilon}^\prime)$ if
$\underline{\epsilon}=\underline{\epsilon}^\prime$ and
 $|\underbar{u}><\underbar{v}|=|\underbar{p}>< \underbar{w}|\in
 \B(\mathbf h^{\otimes n}).$ Expanding the vectors in term of
orthonormal basis $\{e_{\underbar{j}}=e_{j_1}\otimes \cdots
\otimes  e_{j_n}: \underbar{j} =(j_1, \cdots, j_n) , j_1, \cdots,
j_n \ge 1\}$ from $\mathcal D(G)$, the identity
$|\underbar{u}><\underbar{v}|=|\underbar{p}>< \underbar{w}|$ is
equivalent to
 ${\underbar{u}}_{\underbar
{j}}{\overline{\underbar{v}}}_{\underbar {k}}
    ={\underbar{p}}_{\underbar {j}}
{\overline{\underbar{w}}}_{\underbar {k}}
    $ for each
multi-indices $\underbar {j},\underbar {k}$ which gives, $
 (\underbar{u},\underbar{v},\underline{\epsilon})\sim (\underbar{p},\underbar{w},\underline{\epsilon}^\prime)
\Leftrightarrow
A^{(\underline{\epsilon})}(\underbar{u},\underbar{v})
=A^{(\underline{\epsilon}^\prime)}(\underbar{p},\underbar{w})$ for
all bounded operator $A$  and make  $``\sim"$ a well defined
equivalence relation.
 Now consider the
algebra  $M$ generated  by  $M_0/\sim$ with  multiplication
structure given by
$(\underbar{u},\underbar{v},\underline{\epsilon}).
(\underbar{p},\underbar{w},\underline{\epsilon}^\prime)=
(\underbar{u} \otimes \underbar{p},\underbar{v} \otimes
\underbar{w},\underline{\epsilon}\oplus \underline{\epsilon}^\prime)
.$  We
 define a scalar
valued map $K$ on $M\times M $ by setting, for $(\underbar{u}
,\underbar{v},\underline{\epsilon}),
(\underbar{p},\underbar{w},\underline{\epsilon}^\prime) \in M_0,$
\[K\left( (\underbar{u},\underbar{v},\underline{\epsilon}), (\underbar{p},
\underbar{w},\underline{\epsilon}^\prime)   \right
):=\lim_{t\rightarrow 0}\frac{1}{t} \langle
(U_t^{(\underline{\epsilon})}-1)(\underbar{u}, \underbar{v})
\Omega,~ (U_t^{\underline{\epsilon}^\prime}-1)(\underbar{p},
\underbar{w})~\Omega \rangle,~\mbox{if it exists.}~\]

\bpro \label{noise} If  $\mathcal N_0\subseteq \mathcal D(\mathcal
L)$ then we
have the following.\\
{\bf(i)} The map $K$ is a well defined positive definite
kernel on $M.$\\
{\bf(ii)} Up to unitary equivalence there exists a unique separable
Hilbert space $\mathbf k,$ an  embedding $\eta: M\rightarrow \mathbf
k$ and a representation $\pi$ of  $M,~ \pi:M\rightarrow \mathcal B(
\mathbf k)$ such that
 \be \label{eta-dense11}
\{\eta (\underbar{u},\underbar{v},\underline{\epsilon}):
(\underbar{u},
 \underbar{v},\underline{\epsilon})\in M_0\} ~\mbox{ is total in}~
 \mathbf k,\ee \be \label{eta-kelnel-11} \langle \eta
(\underbar{u},\underbar{v},\underline{\epsilon}),
\eta(\underbar{p},\underbar{w},\underline{\epsilon}^\prime)  \rangle
= K\left( (\underbar{u},\underbar{v},\underline{\epsilon}),
(\underbar{p},\underbar{w},\underline{\epsilon}^\prime)   \right
)\ee and
 \be \label{repn11} \pi
(\underbar{u},\underbar{v},\underline{\epsilon})\eta
(\underbar{p},\underbar{w},\underline{\epsilon}^\prime)=\eta
(\underbar{u}\otimes \underbar{p},\underbar{v}\otimes
\underbar{w},\underline{\epsilon}\oplus
\underline{\epsilon}^\prime)- \langle
\underbar{p},\underbar{w}\rangle \eta
(\underbar{u},\underbar{v},\underline{\epsilon}).\ee
 {\bf(iii)} For
any $(\underbar{u},\underbar{v},\underline{\epsilon})\in M_0,$
$\underbar{u}=\otimes_{i=1}^n u_i,\underbar{v}=\otimes_{i=1}^n v_i$
and $ \underline{\epsilon}=(\epsilon_1,\cdots, \epsilon_n)$
 \be \label{eta-n11}
\eta(\underbar{u},\underbar{v},\underline{\epsilon})= \sum_{ i=1}^n
 \prod_{k\ne i} \langle u_k,v_k \rangle \eta(u_i,v_i, \epsilon_i)
 \ee
{\bf(iv)}  $\eta(u,v,1)=-\eta(u,v,0),~\forall u \in \mathbf h,v\in
\mathcal D(G).$\\
 \noindent {\bf(v)} Writing $\eta(u,v)$ for the vector $\eta(u,v,0)\in
\mathbf k,$ \be \overline{Span}\{ \eta(u,v): u \in \mathbf h ,v\in
\mathcal D(G)\}=\mathbf k.\ee

 \epro
\begin{proof} {\bf(i)}
First note that for any
$(\underbar{u},\underbar{v},\underline{\epsilon})\in M_0,$
$\underbar{u}=\otimes_{i=1}^n u_i,\underbar{v}=\otimes_{i=1}^n
v_i,\\
\underline{\epsilon}=(\epsilon_1,\cdots, \epsilon_n)$ we can write
 \bea \label{Ut-1-decompose11}
&&(U_t^{(\underline{\epsilon})}-1)(\underbar{u},\underbar{v})=
\prod_{i=1}^n U_t^{(\epsilon_i)}(u_i,v_i)-\prod_{i=1}^n \langle u_i,v_i \rangle {\nonumber}\\
&&=\sum_{1\le i\le n} (U_t-1)^{(\epsilon_i)}(u_i,v_i)\prod_{j\ne i} \langle u_j,v_j \rangle{\nonumber}\\
&&+\sum_{2\le l\le n} ~~\sum_{1\le i_1<\cdots < i_m\le n}
\prod_{k=1}^l(U_t-1)^{\epsilon_{i_k}}(u_{i_k},v_{i_k}) \prod_{j\ne
i_k} \langle u_j,v_j \rangle.\eea

\noindent Now by Lemma \ref{4Ut-111}, for elements
$(\underbar{u},\underbar{v},\underline{\epsilon}),
(\underbar{p},\underbar{w},\underline{\epsilon}^\prime) \in M_0,$
$\underline{\epsilon}\in \mathbb Z_2^m$  and
$\underline{\epsilon}^\prime\in \mathbb Z_2^n,$
 we have
\bea \label{KK} && K\left(
(\underbar{u},\underbar{v},\underline{\epsilon}), (\underbar{p},
\underbar{w},\underline{\epsilon}^\prime)   \right
)=\lim_{t\rightarrow 0}\frac{1}{t} \langle
(U_t^{(\underline{\epsilon})}-1)(\underbar{u}, \underbar{v})
\Omega,~ (U_t^{\underline{\epsilon}^\prime}-1)(\underbar{p},
\underbar{w})~\Omega \rangle\\
&&= \sum_{1\le i\le m,~1\le j\le n}
 \prod_{k\ne i} \overline{\langle u_k,v_k \rangle}
 \prod_{l\ne j} \langle p_l,w_l \rangle
 \lim_{t\rightarrow 0}\frac{1}{t}
\langle (U_t-1)^{(\epsilon_i)}(u_i,v_i) ~\Omega,
(U_t-1)^{\epsilon_j^\prime}(p_j,w_j) ~\Omega\rangle \nonumber. \eea

\noindent We note  that \bean && \langle (U_t-1)(u,v)~\Omega,
(U_t-1)(p,w) ~\Omega\rangle\\
&&=\langle U_t(u,v)\Omega,~ U_t(p,w) ~\Omega \rangle-
\overline{\langle u,v\rangle} \langle
p,w\rangle \\
&&- \overline{\langle  u,v\rangle} \langle    \Omega,~
 [(U_t-1)(p,w)]~
 \Omega \rangle\\
&&-\overline{\langle \Omega,~ [(U_t-1)(u,v)]
 \Omega \rangle} \langle  p,w\rangle\\
 &&=\langle p, (Z_t-1)(|w><v|) u \rangle - \overline{\langle  u,v\rangle} \langle
 p,
 [(T_t-1) w \rangle-\overline{\langle u, (T_t-1)v \rangle} \langle  p,w\rangle. \eean Thus existence of the  limits on the right hand
side of (\ref{KK}) follows from the identity (\ref{UtUt*inner})
since the semigroups $T_t$ on $\mathbf h$  and $Z_t$ on
$\B_1(\mathbf h)$ are strongly continuous  and  $|w><v|$ is in
$\mathcal D (\mathcal L).$ Hence $K$ is well defined on $M_0.$ Now
extend this to the algebra $M$ sesqui-linearly.  In particular we
have
 \bea \label{kernel11}
&& K( (u,v,\epsilon ), (p,w,\epsilon^\prime) ) {\nonumber} \\
&&= (-1)^{\epsilon+\epsilon^\prime} \lim_{t\rightarrow 0}\{ \langle
p,\frac{Z_t-1}{t}(|w><v|)u\rangle
 - \overline{\langle u, v \rangle }~\langle  p,\frac{T_t-1}{t} w \rangle
-\overline{\langle u, \frac{T_t-1}{t}v \rangle}~\langle   p,w
\rangle \} {\nonumber}\\
 &&= (-1)^{\epsilon+\epsilon^\prime} \{ \langle p,\mathcal L (|w><v|) u \rangle-\overline{\langle u,v
\rangle} \langle p,G~w  \rangle- \overline{\langle u, G ~v
\rangle}\langle p,w  \rangle\}. \eea

\noindent  Positive definiteness  is obvious as  in \cite{SSS}.\\
 \noindent {\bf (ii)} The  Kolmogorov's
construction \cite{krp} to the pair $(M,K)$  provides the separable
Hilbert space $\mathbf k$ as span closure  of  $\{\eta
(\underbar{u},\underbar{v},\underline{\epsilon}):(\underbar{u},\underbar{v},\underline{\epsilon})
\in M_0 \}.$  Now defining   $\pi$ by (\ref{repn11}) we obtain a
 representation of the algebra $M$ in $\mathbf k$ (proof goes similarly as in Lemma 7.1
\cite{SSS}.\\
 \noindent {\bf(iii)} For any
$(\underbar{p},\underbar{w},\underline{\epsilon}^\prime)\in M_0,$
 by (\ref{Ut-1-decompose11})  and Lemma \ref{4Ut-111},~ we have
\bean && \langle
\eta(\underbar{u},\underbar{v},\underline{\epsilon}),
\eta(\underbar{p}, \underbar{w},\underline{\epsilon}^\prime)
\rangle=K\left( (\underbar{u},\underbar{v},\underline{\epsilon}),
(\underbar{p}, \underbar{w},\underline{\epsilon}^\prime)
\right)\\
&&=\lim_{t\rightarrow 0}\frac{1}{t} \langle
(U_t^{(\underline{\epsilon})}-1)(\underbar{u}, \underbar{v})
\Omega,~ (U_t^{\underline{\epsilon}^\prime}-1)(\underbar{p},
\underbar{w})~\Omega \rangle\\
&&= \sum_{i=1}^n
 \prod_{k\ne i} \overline{\langle u_k,v_k \rangle}
 \lim_{t\rightarrow 0}\frac{1}{t}
\langle (U_t-1)^{(\epsilon_i)}(u_i,v_i)
~\Omega,(U_t^{\underline{\epsilon}^\prime}-1)(\underbar{p},
\underbar{w}) ~\Omega\rangle\\
&&= \sum_{i=1}^n
 \prod_{k\ne i} \overline{\langle u_k,v_k \rangle}
\langle \eta(u_i,v_i,\epsilon_i),\eta(\underbar{p},
 \underbar{w},\underline{\epsilon}^\prime)\rangle.
 \eean
Since $\{\eta(\underbar{p},
\underbar{w},\underline{\epsilon}^\prime):(\underbar{p},
\underbar{w},\underline{\epsilon}^\prime)\in M_0\}$ is a total
subset of $\mathbf k,$   (\ref{eta-n11}) follows.\\

\noindent {\bf(iv)} By (\ref{UtUt*inner}) we have \[\langle
\eta(u,v,1),\eta(\underbar{p},
\underbar{w},\underline{\epsilon}^\prime)\rangle=\langle-
\eta(u,v,0),\eta(\underbar{p},
\underbar{w},\underline{\epsilon}^\prime)\rangle.\] Since
$\{\eta(\underbar{p},
\underbar{w},\underline{\epsilon}^\prime):(\underbar{p},
\underbar{w},\underline{\epsilon}^\prime)\in M_0\}$ is a total
subset of $\mathbf k, \eta(u,v,1)=-\eta(u,v,0).$ \noindent {\bf(v)}
It follows immediately  from parts {\bf(iii)}  and {\bf(iv)}.
\end{proof}
\brmrk \label{trivial-rep} The  representation $\pi$ of $M$ in
$\mathbf k$ is trivial \be  \label{trivial} \pi
(\underbar{u},\underbar{v},\underline{\epsilon})
 \eta(\underbar{p},\underbar{w},\underline{\epsilon}^\prime)=
 \langle \underbar{u},\underbar{v}\rangle
 \eta(\underbar{p},\underbar{w},\underline{\epsilon}^\prime). \ee
 If we redefine $\mathcal M$  to be generated by $\underbar{u},\underbar{v} \in \mathcal D(G)^{\otimes n},$
 then  $\mathcal M$ can be a $*$-algebra with  involution:  $(\underbar{u},\underbar{v},\underline{\epsilon})^*=
 (\underleftarrow{u},\underleftarrow{v},\underline{\epsilon}^*)$  (for notations see \cite{SSS} )
  and it is obvious that  $\pi$  given by (\ref{trivial}) is indeed a $*$-representation.\ermrk

 \noindent In the sequel, we   fix  an orthonormal basis  $\{E_j:j\ge 1\}$ of  $\mathbf k.$
 \blema \label{H,Lj11} Under the hypothesis  of Proposition \ref{noise}
we have the followings.
\begin{description}
\item[(i)] There exists a unique family of operators $\{L_j:j\ge 1\}$ in $\mathbf
h$ with  $\mathcal D(L_j)\supseteq \mathcal D(G)$
 such that $\langle u,L_j v \rangle=\eta_j(u,v):=\langle E_j,
\eta(u,v) \rangle, \forall u\in \mathbf h, v\in \mathcal D(G)$ and \\
$ \sum_{j\ge 1 } \|L_j v\|^2= -2~ Re~\langle v, G~v \rangle,~
\forall ~v\in \mathcal D(G).$
\item[(ii)]
The family of operators $\{L_j:j\ge 1\}$  satisfies
 $\sum_{j\ge 1}  \langle u, c_j L_j v\rangle =0, \forall  u \in \mathbf h , v\in \mathcal D(G)$
 for some $c=(c_j)\in l^2(\mathbb N)$
  implies $c=0.$
\item[(iii)]
The generator  $\mathcal L$ of strongly continuous semigroup $Z_t$
satisfies
 \be \label{LZ11}
  \langle p,\mathcal L (|w><v|) u \rangle=
\langle p, |Gw><v| ~u\rangle + \langle p,|w><Gv| ~u\rangle
+\sum_{j\ge 1} \langle p, |L_j w>< L_j v|~u\rangle,\ee for all $u,p
\in \mathbf h  $  and $ v,w\in \mathcal D(G).$ Furthermore, the
family of operators $G, L_j:j\ge 1$  satisfies \be \label{fel}
\langle v, Gw \rangle + \langle Gv, w \rangle +\sum_{j\ge 1} \langle
L_j v, L_j w\rangle = 0, \ee for all $ v,w \in \mathcal D(G).$
\end{description}
\elema
\begin{proof}
{\bf(i)} By the identity (\ref{kernel11}),   for any $u \in \mathbf
h , v\in \mathcal D(G)$ \bea \label{etanorm}
&&\|\eta(u,v)\|^2   {\nonumber}\\
&&= \langle u,\mathcal L (|v><v|) u \rangle-\overline{\langle u,v
\rangle} \langle u,G~v  \rangle- \overline{\langle u, G ~v
\rangle}\langle u,v  \rangle\\
 &&\le \{ \|\mathcal L (|v><v|)\|_1 +2 \| G ~v\| ~\|v\| \}~\|u\|^2 {\nonumber}. \eea

\noindent Thus the linear map $ \mathbf h  \ni u\mapsto \eta(u,v)
\in \mathbf k$ is  a bounded linear map.  Hence by Riesz's
representation theorem, there exists  unique linear operator $L$
from $\mathcal D(G) $  to $\mathbf h\otimes \mathbf k$ such that
$\langle  \langle
 u,Lv \rangle  \rangle   = \eta (u,v)$ where the vector  $\langle  \langle
 u,Lv \rangle  \rangle\in \mathbf k$ is defined as in (\ref{partinn}).
 Equivalently, there exists a unique family
of linear operator $\{L_j:j\ge 1\} $ from $\mathcal D(G)$  to
$\mathbf h$ such that $Lu=\sum_{j\ge 1} L_j u \otimes E_j$ and
$\langle u,L_j v \rangle=\eta_j(u,v).$ Now, for any  $v\in \mathcal
D(G)$ \bean && \|Lv\|^2= \sum_{j} \|L_j v \|^2=\sum_{j,k}
|\eta_j(e_k,v)|^2
=\sum_{k} \|\eta(e_k,v)\|^2\\
&&=\sum_{k}\left [\langle e_k, \mathcal L(|v><v|)
e_k\rangle-\overline{\langle e_k, v\rangle} \langle e_k, G~ v\rangle
- \overline{\langle e_k, G~v  \rangle} \langle e_k, v \rangle ~~\right ]\\
&&=Tr \mathcal L(|v><v|)  -\langle v, G~v \rangle-\overline{\langle
v, G~v \rangle}.\eean

\noindent Since $Z_t$ is trace preserving (\ref{tr-zero11}) and
$|v><v|\in \mathcal D(\mathcal L) $ by hypothesis it follows that
\[Tr \mathcal L(|v<v|)=0\] and therefore \be \label{lj*lj11}
\|Lv\|^2= \sum_{j} \|L_j v \|^2 = -\langle v, G~v
\rangle-\overline{\langle v, G~v \rangle}= -2 Re \langle v, G~v
\rangle. \ee  Note that the term on right hand side  is positive
since   $G$ is the generator of a contractive semigroup.

\noindent{\bf(ii)} For some $c=(c_j)\in l^2(\mathbb N)$  let
$\langle u, \sum_{j\ge 1} c_j L_j v \rangle =0,~ \forall ~ u  \in
\mathbf h ,v\in \mathcal D(G).$ We have
\[0=\langle u, \sum_{j\ge 1} c_j L_j v \rangle = \sum_{j\ge 1} c_j
\langle u, L_j v \rangle =\langle
 \sum_{j\ge 1} \overline{c}_j E_j,\eta (u,v) \rangle.\]
Since $\overline{Span}\{ \eta(u,v): u  \in \mathbf h ,v\in \mathcal
D(G)\}=\mathbf k ,$ it follows that $
 \sum_{j\ge 1} \overline{c}_j E_j=0\in \mathbf k$  and hence $c_j=0,~\forall
 j.$\\
\noindent{\bf(iii)} By  part {\bf(i)}  and identity
(\ref{kernel11}), for any $u,p \in \mathbf h $ and $ v,w\in \mathcal
D(G)$ we have
 \bean
&&   \sum_{j\ge 1} \overline{ \langle u, L_j v \rangle } \langle p,
L_j w \rangle=\langle \eta(u,v), \eta(p,w) \rangle
 \\
 &&= \langle p,\mathcal L (|w><v|) u \rangle-\overline{\langle u,v
\rangle} \langle p,G~w  \rangle- \overline{\langle u, G ~v
\rangle}\langle p,w  \rangle. \eean Thus \bean
&&\langle p, \mathcal L (|w><v|)~u\rangle \\
&&= \langle p, |Gw><v| ~u\rangle + \langle p,|w><Gv| ~u\rangle
+\sum_{j\ge 1} \langle p, |L_j w>< L_j v|~u\rangle. \eean

\noindent  Since, for any $v,w\in \mathcal D(G),$ by identity
(\ref{tr-zero11}),  $ Tr [\mathcal L (|w><v|)]= 0,$ from the above
identity
%(\ref{LZ11})
we get \be \label{fel0}\langle v, Gw \rangle + \langle Gv, w \rangle
+\sum_{j\ge 1} \langle L_j v, L_j w\rangle= 0.\ee

\end{proof}
\brmrk If there exists a positive self adjoint operator  $A$  such
that $\langle v, Av \rangle=-2Re \langle v,Gv \rangle,~ \forall
v\in\mathcal D(G),$ then $\|Lv\|^2= \sum_{j} \|L_j v \|^2= \langle
v, Av \rangle=\|A^\half v \|^2, ~ \forall v\in \mathcal
D(G)\subseteq \mathcal D(A) \subseteq \mathcal D(A^\half)$ and hence
$L$ will be closable. Closability of $(L, \mathcal D(G))$ can be
seen as follows. Suppose $\{v_n\}\subseteq \mathcal D(G) $ converges
to  $0$ and $\{Lv_n\}$ is convergent. Since $\|L(v_n-v_m) \|=
\|A^\half (v_n-v_m) \|,$ convergence  of $\{Lv_n\}$ implies
$\{A^\half v_n\}$ is Cauchy, so convergent in $\mathbf h.$ As
$A^\half$ is a closed operator we get that $A^\half v_n$ converges
to $0$ which implies $Lv_n$ converges to $0.$

This can happen e.g.  when
 $\{T_t\}$  is a   holomorphic semigroup of contractions.\ermrk

 \brmrk \label{D-DG}
 If we replace  $\mathcal D(G)$  by  any dense  subset $\mathcal D\subseteq \mathcal D(G),$ such that  $|u><v|\in \mathcal D(\mathcal
 L)$  for all $u,v\in \mathcal D,$  then  above Proposition \ref{noise} and Lemma \ref{H,Lj11}
 hold with the tensor algebra $\mathcal M$ modified so as to be
 generated by   $  (\otimes_{i=1}^n u_i,\otimes_{i=1}^n v_i ) : u_i
   \in \mathbf h$  and $v_i\in \mathcal D.$
 \ermrk
\section{Hudson-Parthasarathy (HP)  Flows  and  Equivalence}

\noindent In order to set up the  Hudson-Parthasarathy (HP) equation
and proceed further
 we shall work under the following  extra  assumption.\\
\noindent {\bf Assumption E:} There exists a dense set
 $\mathcal D \subseteq \mathcal D(G) \cap \mathcal D(G^*)$ such that
 $\mathcal D$ is a core
of $G$ in $\mathbf h$   and
\begin{description}
\item [E1.] $\mathcal D\subseteq \mathcal D(L_j^*)$ for every  $j\ge 1,$
\item [E2.]
$\mathcal N=Span\{|u><v| : u,v \in \mathcal D\}$ is  core for
  the generator $\mathcal L$  and $\mathcal L^\prime$ of the
semigroup $Z_t$   and $F_t$ on $\B_1(\mathbf h)$  respectively,
\item [E3.] $L_j$  maps  $\mathcal D$  into itself  and
for any $v\in \mathcal D,  \sum_{j\ge 1} \| G L_j v\|^2 <\infty.$
\end{description}
Since  $\mathcal D$  is dense in $\mathbf h$ one can see,  by a
simple approximation argument,  that $\mathcal N $  is dense in
$\B_1(\mathbf h).$  Recall from the Remark \ref{D-DG} that under the
 assumption {\bf E2}, replacing $\mathcal D(G) $  by  the core
$\mathcal D$ in Proposition \ref{noise} and Lemma \ref{H,Lj11}, we
get a separable Hilbert space $\mathbf k$ generated  by $\{
\eta(u,v): u \in \mathbf h ,v\in \mathcal D\} $ and linear operators
$ \{L_j:j\ge 1\}$ defined on $\mathcal D.$

\brmrk The  assumption {\bf E1}  is  needed for setting up an  HP
equation with coefficients  $G$ and $L_j: j\ge 1,$   assumption {\bf
E2} is to assure the   existence of unique  unitary   HP  flow. The
 assumption {\bf E3} will be necessary for proving the minimality of
the associated HP flow which will be needed  to establish  unitary
equivalence  of   the  HP  flow  and   unitary process $U_t,$ we
started with. \ermrk

\noindent Now let us state the main result of this article.
\bthm \label{mainthm} Assume {\bf A,B, C, D and E}. Then   we have the following.\\
{\bf (i)} The  HP  equation  \be \label{hpeqn,st11} V_t=1_{\mathbf h
\otimes \Gamma}+\sum_{\mu,\nu\ge 0}\int_0^t V_r L_\nu^\mu
\Lambda_\mu^\nu(dr)\ee  on $\mathcal D \otimes \mathcal
E(L^2(\mathbb R_+,\mathbf k))$  with coefficients $L_\nu^\mu$ given
by
\begin{equation} \label{hpcoefi11}
L_\nu^\mu=\left\{ \begin{array} {lll}
 & G & \mbox{for}\
(\mu,\nu)=(0,0)\\
&L_j &  \mbox{for}\ (\mu,\nu)=(j,0)
\\
& -L_k^*&
\mbox{for}\ (\mu,\nu)=(0,k)\\
& 0 & \mbox{for}\ (\mu,\nu)=(j,k)
 \end{array}
 \right.
  \end{equation}
admit a unique unitary solution $V_t.$

\noindent
 {\bf (ii)} There exists a unitary isomorphism $\widetilde{\Xi}:\mathbf h \otimes \mathcal H
\rightarrow \mathbf h \otimes \Gamma$  such that \be \label{U=V} U_t
=\widetilde{\Xi}^*~ V_t~ \widetilde{\Xi},~\forall~~ t\ge 0.\ee \ethm

 \noindent Here we shall sketch the prove of   part {\bf (i)} of
the Theorem and postponed the proof of  {\bf (ii)} to  next two sub
sections. In order to prove the part  {\bf (i)}  we need  the
following  two Lemmas. For $\lambda >0,$ we define the Feller set  $
\beta_\lambda\subseteq \B(\mathbf h)$ by\\ $ \{ x\ge 0 :\langle v, x
L_0^0 w\rangle + \langle  L_0^0 v,xw \rangle +\sum_{j\ge 1} \langle
L_0^j v,x L_0^j w\rangle =\langle v, xGw\rangle + \langle Gv,xw
\rangle +\sum_{j\ge 1} \langle L_j v,x L_j w\rangle=\lambda \langle
v,x w\rangle, \forall v,w\in \mathcal D\}.$ Similarly we  define the
Feller set $ \widetilde{\beta_\lambda}  $  for  coefficients
$\widetilde {L}_\nu^\mu\equiv (L_\mu^\nu)^*.$

 \blema \label{LemaFel} Under the  assumption E2,  the Feller Condition:  $
 \beta_\lambda=\{0\}$    as well as  $
 \widetilde{\beta_\lambda}=\{0\}$ for some $\lambda >0$ hold
 .
 \elema
 \begin{proof} For any $x\ge 0$ in $\B(\mathbf h), v,w\in \mathcal D$
we have
 \bean &&\sum_{j\ge 1} \langle L_j v,x L_j w\rangle= \langle L v,x L w\rangle
 =\langle x^\half  L v,x^\half  L w\rangle=\sum_{m
\ge 1} \langle  L v,(|x^\half e_m ><  x^\half  e_m| \otimes
1_{\mathbf k})L w
 \rangle\\
&&=\sum_{m\ge 1} \langle ~\langle \langle x^\half e_m,  L  v \rangle
\rangle ,  \langle \langle x^\half e_m, L w \rangle \rangle
~\rangle= \sum_{m\ge 1} \langle \eta( x^\half e_m,v),\eta( x^\half
e_m,w)\rangle.\eean Now by (\ref{kernel11})
 \bea  \label{innFel} &&\sum_{j\ge 1} \langle L_j v,x L_j
w\rangle= \sum_{m\ge 1} \langle \eta(  x^\half e_m,v),\eta(
 x^\half e_m,w)\rangle\\
&& = \sum_{m\ge 1}\{ \langle  x^\half e_m, \mathcal L(|w><v|)
x^\half e_m \rangle- \overline{\langle  x^\half e_m, G v  \rangle}
\langle x^\half e_m ,w \rangle- \overline{\langle  x^\half e_m, v
\rangle} \langle  x^\half e_m , G w \rangle  \} {\nonumber}\\
&& = Tr [x  \mathcal L(|w><v|)]- \langle  v,xG w\rangle -\langle G
v, xw  \rangle  {\nonumber}.\eea Thus  \be \label{LG} \langle v,
xGw\rangle + \langle Gv,xw \rangle +\sum_{j\ge 1} \langle L_j v,x
L_j w\rangle=Tr [x \mathcal L(|w><v|)]\ee  and for any  $x\in
\beta_\lambda,$ \be \label{fe} Tr [x \mathcal L(|w><v|)] = \lambda
\langle v, xw\rangle= \lambda~ Tr (x |w><v|), \forall v,w\in
\mathcal D.\ee By assumption  {\bf E2}  the subspace  $\mathcal
N=Span\{|w><v|:v,w\in \mathcal D\}$ is a core for $ \mathcal L$ and
hence  the identity (\ref{fe}) extends to $ Tr [x \mathcal L(\rho)]
=\lambda~ tr (x\rho), \forall \rho\in \mathcal D(\mathcal L).$ It is
also clear that for  $x\in \beta_\lambda$ the scalar map $ \phi_x:
 \mathcal D(\mathcal L) \ni \rho \mapsto Tr [x \mathcal L(\rho)]=\lambda~  Tr (x \rho) $ extends to a bounded linear
functional   on $\B_1(\mathbf h).$   Hence $x$ is in the domain of
$\mathcal L^*$ and we get \bean &&Tr [(|w><v|)(\mathcal
L^*-\lambda)x]=0\\
&&\Rightarrow \langle v, (\mathcal L^*-\lambda)x w \rangle=0\\
&&\Rightarrow (\mathcal L^*-\lambda)x=0.\eean Since $\mathcal L^*$
is the generator of a $C_0$-semigroup  $\{Z_t^*\}$ of contraction
maps on $\B(\mathbf h),$ for $\lambda >0, ~ \mathcal L^*-\lambda $
is invertible  and hence $x=0.$

\noindent  To prove $\widetilde{\beta_\lambda}=\{0\}$  let us
consider the following.  By identity (\ref{UtUt*10}) for  vectors
$u,p\in \mathbf h$ and $v,w\in \mathcal D$ \bean && \langle \eta(
u,v),\eta( p,w)\rangle \\
&& = \lim_{t\rightarrow 0}\frac{1}{t} \langle (U_t-1)(u,v)
~\Omega, (U_t-1)(p,w) ~\Omega\rangle\\
&&= \lim_{t\rightarrow 0}\frac{1}{t} \langle (U_t^*-1)(u,v) ~\Omega,
(U_t^*-1)(p,w) ~\Omega\rangle\\
&&= \lim_{t\rightarrow 0} \frac{1}{t}\{\langle U_t^*(u,v)\Omega,~
U_t^*(p,w) ~\Omega \rangle- \overline{\langle u,v\rangle} \langle
p,w\rangle \\
&&- \overline{\langle  u,v\rangle} \langle    \Omega,~
 [(U_t^*-1)(p,w)]~
 \Omega \rangle\\
&&-\overline{\langle \Omega,~ [(U_t^*-1)(u,v)]
 \Omega \rangle} \langle  p,w\rangle\}\\
 &&=  \lim_{t\rightarrow 0} \frac{1}{t} \{\langle p, (F_t-1)(|w><v|) u \rangle - \overline{\langle  u,v\rangle} \langle
 p,
 (T_t^*-1) w \rangle-\overline{\langle u, (T_t^*-1)v \rangle} \langle
 p,w\rangle\}.\eean
 Since by E2,~ $v,w\in \mathcal D\subseteq  \mathcal D(G^*)$  and  $|w><v| \in  \mathcal D(\mathcal
 L^\prime),$ we get that
  \be \label{eta*inn}\langle \eta(
u,v),\eta( p,w)\rangle=\langle p, \mathcal L^\prime(|w><v|) u
\rangle - \overline{\langle  u,v\rangle} \langle
 p,
 G^* w \rangle-\overline{\langle u, G^*v \rangle} \langle
 p,w\rangle.\ee
Thus by (\ref{innFel})  and  (\ref{eta*inn}) we have
 \bean &&\sum_{j\ge 1} \langle L_j v,x L_j
w\rangle= \sum_{m\ge 1} \langle \eta(  x^\half e_m,v),\eta(
 x^\half e_m,w)\rangle\\
&& = \sum_{m\ge 1}\{ \langle  x^\half e_m, \mathcal L^\prime
 (|w><v|) x^\half e_m \rangle- \overline{\langle  x^\half
e_m, v \rangle} \langle  x^\half e_m , G^* w \rangle -
\overline{\langle  x^\half e_m, G^* v
\rangle} \langle  x^\half e_m ,w \rangle \}\\
&& = Tr [x \mathcal L^\prime (|w><v|)] -\langle G^* v, xw \rangle-
\langle  v,xG^* w\rangle.\eean Thus  \be \label{LG*} \langle v,
xG^*w\rangle + \langle G^*v,xw \rangle +\sum_{j\ge 1} \langle L_j
v,x L_j w\rangle=Tr [x \mathcal L^\prime (|w><v|)]\ee  and for any
$x\in \widetilde{\beta_\lambda},$ \be \label{fe*} Tr [x \mathcal
L^\prime(|w><v|)] = \lambda \langle v, xw\rangle= \lambda~ Tr
(x|w><v|), \forall v,w\in \mathcal D.\ee

Since  the subspace $\mathcal N=Span\{|w><v|:v,w\in \mathcal D\}$ is
a core for $ \mathcal L^\prime$ by    assumption {\bf E2}, a similar
argument  as above will give  that
$\widetilde{\beta_\lambda}=\{0\}.$

 \end{proof}

\brmrk By  (\ref{LG})  and (\ref{LG*}) formally  $(\mathcal
L^\prime-\mathcal L )\rho= [G^*-G, \rho], \forall \rho \in \mathcal
N. $ Denoting  the imaginary part of  $G$ by $H$ consider the
derivation $ \delta_H (\rho)=-2 ~i~[H, \rho].$ If $\delta_H$ is
bounded then the hypothesis that the subspace
 $\mathcal N$ is a core for  $\mathcal L$ implies that it
is a core for $\mathcal L^\prime$ and  no extra assumption is
needed. \ermrk

\brmrk If
 $\{T_t\}$  is a   holomorphic semigroup of contractions then  the hypotheses   on  domains of $G^*$
  and  $\mathcal L^\prime$   will hold automatically.

\ermrk

\blema Assume  the hypotheses  {\bf E1  }  and   {\bf E2  }.  For
$n\ge 1,$ setting $L_j(n)=n~L_j~(n 1_{\mathbf h}-G)^{-1} $ and
$G(n)=n^2 (n 1_{\mathbf h}-G^*)^{-1}  G (n 1_{\mathbf
h}-G)^{-1},$ we have. \\
(i) The operators $L_j(n), G(n) \in \B(\mathbf h)$ and
$\sum_{j} \|L_j(n) v \|^2=-2~Re  \langle  v, G(n) v \rangle.$\\
(ii) For  $v\in \mathcal D,\lim_{n \rightarrow \infty} L_j(n)v=L_j
v,~\lim_{n \rightarrow \infty} L_j(n)^*v=L_j^* v$  and
\\$~~~~~~\lim_{n \rightarrow \infty} G(n)v=G v.$
 \elema
\begin{proof} (i) For any $v\in \mathbf h,$
\bean && \sum_{j} \|L_j(n) v \|^2=\sum_{j}n^2 \| L_j~(n 1_{\mathbf
h}-G)^{-1} v \|^2 \\
&&=-2 Re~n^2 \langle  (n 1_{\mathbf h}-G)^{-1} v, G (n 1_{\mathbf
h}-G)^{-1} v \rangle\\
&&=-2~Re  \langle  v, G(n) v \rangle.\eean
 \noindent (ii) Since the
sequences of  bounded operators $\{ nL_j  (n 1_{\mathbf h}-G)^{-1}\}
$ and $\{ nL_j  (n 1_{\mathbf h}-G^*)^{-1}\} $ are uniformly norm
bounded  and converge strongly  to identity, the requirements
follows.
\end{proof}
\noindent  {\bf Sketch of the Proof of the part (i) of Theorem
\ref{mainthm} ~:}

\noindent  For each $n\ge 1$ we   consider the  family of operators,
\begin{equation} \label{hpcoefinn}
L_\nu^\mu(n)=\left\{ \begin{array} {lll}
 &G(n)=n^2 (n 1_{\mathbf
h}-G^*)^{-1}G(n 1_{\mathbf h}-G)^{-1} & \mbox{for}\
(\mu,\nu)=(0,0)\\
&L_j(n)=n~L_j~(n 1_{\mathbf h}-G)^{-1}  &  \mbox{for}\
(\mu,\nu)=(j,0)
\\
&-L_k(n)^*&
\mbox{for}\ (\mu,\nu)=(0,k)\\
&0 & \mbox{for}\ (\mu,\nu)=(j,k).
 \end{array}
 \right.
  \end{equation}

 \noindent By hypothesis {\bf E1}, we have  that
$\lim_{n\rightarrow \infty}  L_\nu^\mu(n)v= L_\nu^\mu v, \forall
v\in \mathcal D$  and hence there exist unique contractive  solution
$\{V_t\}$  for the   HP equation (\ref{hpeqn,st11})  (see
\cite{moh,f, gs}). To show that $\{V_t\}$  is a isometric process we
shall  use the Feller condition proved in Lemma  \ref{LemaFel}. By
Proposition 3.1 in \cite{ms} (also see \cite{moh,f}) / Theorem 7.2.3
in \cite{gs} the solution $\{V_t\}$ of HP equation \ref{hpeqn,st11}
is isometric. We shall conclude the unitarity  of the process $V_t$
by employing time reversal operator and  the results in \cite{ms,
gs}. As $V_t$ satisfies  the equation (\ref{hpeqn,st11}), $V_t^*$
satisfies the HP equation on $\mathcal D\otimes \mathcal E(\mathcal
K),$  since $\mathcal D\subseteq  \mathcal D(G^*)$  by E2,
 \be \label{hpeqn,st*}
V_t^*=1_{\mathbf h \otimes \Gamma}+\sum_{\mu,\nu\ge 0}\int_0^t
 (L_\nu^\mu)^* V_r^*
 \Lambda_\nu^\mu(dr).\ee Let us define $\widetilde{V_t}:=  [1_\mathbf h\otimes
\Gamma(R_t)] V_t^*  [1_\mathbf h\otimes  \Gamma(R_t)] ,$  where
$R_t$  is the time reversal operator on $L^2(\mathbb R_+, \mathbf
k):$ \bean && R_t f
(x)= f (t-x) ~ \mbox{if} ~  x\le t\\
&&~~~~~~~~~= f(x) ~~~~~~ \mbox{if} ~x>t \eean and $\Gamma(A)$ denote
the second quantization of operator $A:~
\Gamma(A)\textbf{e}(f)=\textbf{e}(Af).$  Then it can be seen that
the process $\{\widetilde{V_t}\}$ satisfies the HP equation  on
$\mathcal D\otimes \mathcal E(\mathcal K),$
  \be \label{hpeqn,tilde}
\widetilde{V_t}= 1_{\mathbf h \otimes \Gamma}+\sum_{\mu,\nu\ge
0}\int_0^t
 \widetilde{V_r}  \widetilde {L}_\nu^\mu
 \Lambda_\mu^\nu(dr).\ee Since the Feller condition
$ \widetilde{\beta_\lambda}= \{0\}$ for $\widetilde {L}_\nu^\mu$
holds  by Lemma \ref{LemaFel}, the solution $\widetilde{V_t}$ and
hence $V_t^*$ is isometric or equivalently  $V_t$ is co-isometric
and therefore  $V_t$ is a strongly continuous unitary process. \qed

\brmrk
 Using identity (\ref{fel0}) one
 construct the minimal semigroup $\hat{Z_t}$ with generator
$\hat{\mathcal L}$ such that restrictions of $\mathcal L$ and
$\hat{\mathcal L}$  to $\mathcal N$  are same (see \cite{gs,ms,
moh,si}). Therefore, for any $\lambda > 0,$  the closure
$\overline{(\lambda -\hat{\mathcal L} ) \mathcal
N}=\overline{(\lambda -\mathcal L) \mathcal N} =(\lambda -\mathcal
L) \mathcal D(\mathcal L)$ since by hypothesis {\bf E2} the subspace
$\mathcal N$ is a core for $\mathcal L.$ As $\mathcal L$  is the
generator of a $C_0$-semigroup of contractions on $\mathcal
B_1(\mathbf h)$  the subspace $(\lambda -\mathcal L) \mathcal
D(\mathcal L)=\mathcal B_1(\mathbf h) $  and hence
$\overline{(\lambda -\hat{\mathcal L} ) \mathcal N}= \mathcal
B_1(\mathbf h).$ Thus by Theorem 3.2.16 (ii) and (iii) in \cite{gs}
we have that $Tr(\hat{Z_t} \rho)=Tr(\rho),$ i.e the minimal
semigroup $\hat{Z_t}$ is conservative which also implies that the
Feller condition  is satisfied. We also have $\overline{(\lambda
-\hat{\mathcal L} ) \mathcal N}= \mathcal B_1(\mathbf
h)=\overline{(\lambda -\hat{\mathcal L} ) \mathcal D(\hat{\mathcal
L})}$  which implies   $\mathcal N$ is a core for $\hat{\mathcal L}$
as well  and  hence $\mathcal L=\hat{\mathcal L}.$ Thus   $Z_t$  is
the  minimal  semigroup. \ermrk

\noindent  For any $0\le s\le t<\infty,$ we define a unitary
operator $V_{s,t}:= [1_\mathbf h\otimes \Gamma(\theta_s)]  V_{t-s}
[1_\mathbf h\otimes \Gamma(\theta_s^*)],$ where $\theta_s$  is the
right shift operator on $L^2(\mathbb R_+, \mathbf k):$ \bean
&&\theta_s f
(x)= f (x-s) ~ \mbox{if} ~  x\ge s\\
&&~~~~~~~~~= f(x) ~~~~~~ \mbox{if} ~x<s.\eean  The adjoint of $
\theta_s$  is given by $ \theta_s^* f (x)= f (x+s) $ for all  $x\ge
0.$ We shall write the ampliation $1_\mathbf h\otimes A$ of an
operator $A$ by same symbol $A $  when it is clear from the context.
Since the unitary process $V_t$  is the solution of  HP equation
(\ref{hpeqn,st11}) we have
 \bean
&& V_{s,t}= \Gamma(\theta_s)  V_{t-s}
\Gamma(\theta_s^*)\\
&&=1_{\mathbf h \otimes \Gamma}+\sum_{\mu,\nu\ge 0}
\Gamma(\theta_s)\{\int_0^{t-s}
V_r L_\nu^\mu \Lambda_\mu^\nu(dr)\} \Gamma(\theta_s^*)\\
&&=1_{\mathbf h \otimes \Gamma}+\sum_{\mu,\nu\ge 0} \int_0^{t-s}
 \Gamma(\theta_s)  V_r \Gamma(\theta_s^*)  L_\nu^\mu \Gamma(\theta_s) \Lambda_\mu^\nu(dr) \Gamma(\theta_s^*)\\
 &&=1_{\mathbf h \otimes \Gamma}+\sum_{\mu,\nu\ge 0}\int_s^{t}
 \Gamma(\theta_s)  V_{r-s} \Gamma(\theta_s^*) L_\nu^\mu  \Gamma(\theta_s)
 \Lambda_\mu^\nu(d r-s)  \Gamma(\theta_s^*).\eean
Since for any interval  $ \Delta \subseteq \mathbb R_+,
\Gamma(\theta_s) \Lambda_\mu^\nu(\Delta -s)\}
\Gamma(\theta_s^*)=\Lambda_\mu^\nu(\Delta)$
 it follows that  the unitary family $\{V_{s,t}\}$  satisfies the  HP
equation \be \label{hpst} V_{s,t}=1_{\mathbf h \otimes
\Gamma}+\sum_{\mu,\nu\ge 0}\int_s^t V_{s,r} L_\nu^\mu
\Lambda_\mu^\nu(dr)\ee on $\mathcal D\otimes \mathcal E(L^2(\mathbb
R_+, \mathbf k)).$ We note that $V_t=V_{0,t}$  and
$V_{s,s}=1_{\mathbf h \otimes \Gamma}.$

%Also consider \be \label{hpeqn,stnn} V^{(n)}_{s,t}=1_{\mathbf h
%\otimes \Gamma}+\sum_{\mu,\nu\ge 0}\int_s^t V^{(n)}_{s,r}
%L_\nu^\mu(n) \Lambda_\mu^\nu(dr)\ee with coefficients $L_\nu^\mu(n)$
%given by (\ref{hpcoefinn}).

\noindent As for the family of unitary operators $\{ U_{s,t}\}$  on
$\mathbf h \otimes \mathcal H,$ for
$\underline{\epsilon}=(\epsilon_1,\epsilon_2,\cdots, \epsilon_n)\in
\mathbb Z_2^n$ we define $V_{s,t}^{(\underline{\epsilon})}\in
\mathcal B(\mathbf h^{\otimes n}\otimes \Gamma)$  by setting
$V_{s,t}^{(\epsilon)}\in \mathcal B(\mathbf h\otimes \Gamma)$ by
\bean
&& V_{s,t}^{(\epsilon)}=V_{s,t} ~~~\mbox{for}~ \epsilon =0\\
&&~~~~~=V_{s,t}^* ~~~\mbox{for}~ \epsilon =1. \eean The next result
verifies the properties of  assumption {\bf A} for the family
$V_{s,t}$ with $\textbf{e}(0) \in \Gamma$  replacing  $\Omega \in
\mathcal H.$
 \blema \label{Vstbasic1} The family of unitary
operators $\{ V_{s,t}\}$ satisfy
\begin{description}
\item[(i)] For any $0\le r\le  s \le t<\infty, V_{r,t}=V_{r,s} V_{s,t}.$
\item[(ii)] For $[q,r)\cap [s,t)=\emptyset, V_{q,r}(u,v)$ commute with
 $V_{s,t}(p,w)$  and $V_{s,t}(p,w)^*$  for every $u,v,p,w \in \mathbf h.$
\item[(iii)] For any $0\le s \le t<\infty,$\\
$ \langle \textbf{e}(0),   V_{s,t}(u,v) \textbf{e}(0) \rangle =
\langle \textbf{e}(0),   V_{t-s}(u,v) \textbf{e}(0) \rangle =
\langle u,  T_{t-s }v\rangle,~\forall u,v \in \mathbf h.$
\end{description}
\elema
\begin{proof}
{\bf{(i)}} For fixed $0\le r\le  s\le t<\infty,$  we set
$W_{r,t}=V_{r,s} V_{s,t}.$
 Then by (\ref{hpeqn,st11}) we have
\bean
&& W_{r,t}=V_{r,s}+\sum_{\mu,\nu\ge 0}\int_s^t V_{r,s}V_{s,q} L_\nu^\mu  \Lambda_\mu^\nu(dq)\\
&&=W_{r,s}+\sum_{\mu,\nu\ge 0}\int_s^t   W_{r,q} L_\nu^\mu
\Lambda_\mu^\nu(dq). \eean
 Thus  the family of  unitary operators $\{ W_{r,t}\}$ also satisfies the
 HP equation
(\ref{hpst}).  Hence by uniqueness of the solution of this quantum
stochastic differential equation, $W_{r,t}=V_{r,t},
\forall t\ge s$ and the result follows.\\

 \noindent {\bf{(ii)}} For any $0\le s\le t<\infty,~$
$V_{s,t}\in \mathcal B( \mathbf h\otimes \Gamma_{[s,t]}).$
 So for $p,w\in \mathbf h,~V_{s,t}(p,w)\in \mathcal B(\Gamma_{[s,t]})$ and
the statement  follows.\\

  \noindent {\bf{(iii)}} Let us set  a family of contraction operators $\{
  \widetilde{S}_{s,t}\}$ on $\mathbf h$ by
\[ \langle   u, \widetilde{S}_{s,t} v\rangle =\langle u\otimes \textbf{e}(0),
  V_{s,t}v \otimes \textbf{e}(0) \rangle,~\forall u,v \in \mathbf
  h.\] By definition of $V_{s,t},$ we have
$
\langle u\otimes \textbf{e}(0),
  V_{s,t}v \otimes \textbf{e}(0) \rangle
=
  \langle u\otimes \textbf{e}(0),
 \Gamma(\theta_s)  V_{t-s}
\Gamma(\theta_s^*)v \otimes \textbf{e}(0) \rangle
=
\langle u\otimes \textbf{e}(0),
  V_{0, t-s}v \otimes \textbf{e}(0) \rangle
  $ and hence $\widetilde{S}_{s,t}=\widetilde{S}_{0,t-s}.$
 Setting
$\widetilde{S}_{t}:=\widetilde{S}_{0,t}$ the family
$\{\widetilde{S}_{t}:t\ge 0\}$ is a $C_0$-semigroup of contractions
 on $\mathbf h$. Since the unitary process $V_{s,t}$ satisfies the HP equation
 (\ref{hpst}), for  any  $ u,v \in \mathcal D $
\be \label{S_t} \langle   u, \widetilde{S}_{s,t} v\rangle =\langle
u,v \rangle+\int_s^t \langle u,  \widetilde{S}_{s,r}  G v \rangle
dr.\ee Note that $\mathcal D$  is dense core  for  $G$   and
$\widetilde{S}_{s,t}$ is a contractive  family,   so  the  equation
(\ref{S_t})  extend to  $u\in \mathbf h, v\in \mathcal D(G)$  and
hence the family $\{\widetilde{S}_{s,t}\}$
   satisfies the following differential equation
\[  \widetilde{S}_{s,t}=1+\int_s^t \widetilde{S}_{s,r} G dr\]  on the domain $\mathcal D(G).$ Since  $G$ is the generator of
 the  $C_0$-semigroup
$\{T_t\}$  we have
 $\widetilde{S}_{s,t}=\widetilde{S}_{t-s}=T_{t-s}.$ This proves the claim.
\end{proof}

 \noindent Consider the  family of maps $\widetilde{Z}_{s,t}$ defined by
\[\widetilde{Z}_{s,t} \rho = Tr_{\mathcal H}  [V_{s,t} (\rho \otimes
| \textbf{e}(0)>< \textbf{e}(0)|) V_{s,t}^*], ~ \forall \rho \in
\B_1(\mathbf h).\]    As for  $Z_t,$ it can be  seen
   that $\widetilde{Z}_{s,t}$ is a
  contractive family of maps  on $\B_1(\mathbf h)$  and
  in particular, for any $ u,v,p,w \in \mathbf h$
\[\langle   p, \widetilde{Z}_{s,t} (|w><v|) ~u \rangle
 =\langle   V_{s,t} (u,v) \textbf{e}(0),
  V_{s,t} (p,w) \textbf{e}(0) \rangle.\]
  \blema
  The family $\widetilde{Z}_t:=\widetilde{Z}_{0,t} $  is a $C_0$-semigroup
  of contraction  on $\B_1(\mathbf h)$  and
  $\widetilde{Z}_{s,t}=\widetilde{Z}_{t-s}
  =Z_{t-s}.$
  \elema
  \begin{proof}  By (\ref{hpst}) and Ito's formula for $ u,v,p,w \in \mathcal D$
\bean &&\langle   p,[ \widetilde{Z}_{s,t}-1] (|w><v|) ~u \rangle\\
  &&=\langle   V_{s,t} (u,v) \textbf{e}(0),
  V_{s,t} (p,w) \textbf{e}(0) \rangle-\overline{\langle   u,v
   \rangle} \langle  p,w \rangle\\
   &&=\int_s^t \langle   V_{s,\tau} (u,v) \textbf{e}(0),
  V_{s,\tau} (p,G w) \textbf{e}(0) \rangle d\tau
  +\int_s^t \langle   V_{s,\tau} (u,Gv) \textbf{e}(0),
  V_{s,\tau} (p,w) \textbf{e}(0) \rangle d\tau\\
 && +\int_s^t
  \langle   V_{s,\tau} (u,L_jv) \textbf{e}(0),
  V_{s,\tau} (p,L_j w) \textbf{e}(0) \rangle d\tau\\
  &&=\int_s^t \langle   p, \widetilde{Z}_{s,\tau} (|Gw>< v|) ~u \rangle
  d\tau
  +\int_s^t\langle   p, \widetilde{Z}_{s,\tau} (|w><G v|) ~u \rangle  d\tau\\
 &&+\sum_{j\ge 1}\int_s^t\langle   p, \widetilde{Z}_{s,\tau}
  (|L_jw><L_j v|) ~u \rangle
 d\tau.\eean
 Thus
 \be \label{Ztilda} \langle   p,[ \widetilde{Z}_{s,t}-1] (\rho) ~u \rangle
 =\int_s^t \langle   p, \widetilde{Z}_{s,\tau}
 \mathcal L(\rho) ~u \rangle d\tau,
  \ee where $\rho=|w><v|.$
 Since  $\mathcal D$  is dense in  $\mathbf h,~\mathcal N$ is a core for $\mathcal L$ and  $\widetilde{Z}_{s, \tau}
 $ is a contractive family the equation (\ref{Ztilda})  extends to $ u,p\in \mathbf h$  and $\rho \in \mathcal D(\mathcal L).$
Thus the family $\widetilde{Z}_{s,t}$ satisfies the differential
equation
\[\widetilde{Z}_{s,t} (\rho) =\rho+\int_s^t \widetilde{Z}_{s,\tau}
 \mathcal L(\rho)  d\tau,~~\rho\in  \mathcal D(\mathcal L).\]
 Since $\mathcal L$  is the generator  of $C_0$-semigroup  $Z_t,$ it
 follows that
$\widetilde{Z}_{s,t}=\widetilde{Z}_{t-s}
  =Z_{t-s}.$
  \end{proof}

\subsection{Minimality of HP Flows} \noindent In this section  we
shall show the minimality of the HP flow $V_{s,t}$ discussed above
which will  be needed to prove the Theorem  \ref{mainthm} {\bf
(ii)},  i.e, to establish unitary equivalence  of  $U_t$ and $V_t.$
We shall prove here that the subset $ \mathcal S^\prime:=\{ \zeta=
V_{\underbar{s},\underbar{t}}(\underbar{u},\underbar{v})
\textbf{e}(0):=V_{s_1,t_1}(u_1,v_1)\cdots
V_{s_n,t_n}(u_n,v_n)\textbf{e}(0): \underbar{s}=(s_1,s_2, \cdots,
s_n), \underbar{t}=(t_1,t_2, \cdots, t_n)$ $:~ 0 \le s_1\le t_1\le
\cdots \le  s_n\le  t_n< \infty,n\ge 1, \underbar{u}=\otimes_{i=1}^n
u_i \in \mathbf h^{\otimes n},\underbar{v}=\otimes_{i=1}^n v_i\in
\mathcal D^{\otimes n}\}$ is total in  the symmetric Fock space
$\Gamma(L^2(\mathbb R_+,\mathbf k)).$

\noindent
Since  $\mathcal D$ is dense in $\mathbf h$,  by Remark \ref{DS_0} the subset\\
 $ \mathcal S:=\{ \zeta=
U_{\underbar{s},\underbar{t}}(\underbar{u},\underbar{v})
\Omega:=U_{s_1,t_1}(u_1,v_1)\cdots U_{s_n,t_n}(u_n,v_n)\Omega:
\underbar{s}=(s_1,s_2, \cdots, s_n), \underbar{t}=(t_1,t_2, \cdots,
t_n)$ $:~ 0 \le s_1\le t_1\le s_2 \le \ldots\le  s_n\le  t_n<
\infty,n\ge 1, \underbar{u}=\otimes_{i=1}^n u_i \in \mathbf
h^{\otimes n},\underbar{v}=\otimes_{i=1}^n v_i \in \mathcal
D^{\otimes n} \}$ is total in $\mathcal H.$ We also note that  $\{
\eta(u,v): u \in \mathbf h ,v\in  \mathcal D\}$ is total in $\mathbf
k.$

 \blema  Under  the  assumption {\bf E3},  for any  $v\in \mathcal D,  \sum_{i,j\ge 1}
  \|L_i L_j v\|^2 < \infty.$ \elema
  \begin{proof} For any $j\ge 1, L_j v\in \mathcal D$  and by
  Lemma  \ref{H,Lj11}  {\bf (i)},
\[\sum_{i\ge 1}
  \|L_i L_j v\|^2= -\langle L_j v, G L_j v\rangle - \langle G L_j v, L_j v
  \rangle.\] Therefore
  \[ \sum_{i,j \ge 1}
  \|L_i L_j v\|^2= - 2 Re \sum_{j \ge 1}\langle L_j v, G L_j v\rangle
  \le  2  [\sum_{j \ge 1}\| L_j v\|^2]^\half
   [\sum_{j \ge 1}\| G L_j v\|^2]^\half<\infty.\]

  \end{proof}

 \noindent Let $\tau \ge 0$  be fixed. We note that for any $0\le s < t\le \tau,
u\in \mathbf h,v\in \mathcal D$ by HP equation (\ref{hpeqn,st11})
\bea \label{vst-1}
&&\frac{1}{t-s}[V_{s,t}-1] (u,v) \textbf{e}(0) {\nonumber}\\
&& = \frac{1}{t-s}\{\sum_{j\ge 1} \int_s^t V_{s,\lambda} (u,L_jv)
a_j^\dag (d\lambda)+ \int_s^t V_{s,\lambda} (u,Gv)
d\lambda  \}\textbf{e}(0) {\nonumber}\\
&& = \gamma(s,t,u,v)+  \langle u,Gv\rangle
~~\textbf{e}(0)+\zeta(s,t,u,v) +\varsigma(s,t,u,v)  \eea where these
vectors in the Fock space $\Gamma$ are given by
\begin{description}
\item $\gamma(s,t,u,v):=\frac{1}{t-s}\sum_{j\ge 1} \langle u,L_jv\rangle
a_j^\dag([s,t])~~\textbf{e}(0)$
\item $\zeta(s,t,u,v):=\frac{1}{t-s}\sum_{j\ge 1}\int_s^t (V_{s,\lambda}
-1)(u,L_jv) a_j^\dag (d\lambda) ~~\textbf{e}(0)$
\item  $\varsigma(s,t,u,v):=\frac{1}{t-s} \int_s^t (V_{s,\lambda}-1)
(u,Gv) d \lambda ~~\textbf{e}(0).$
\end{description}

\noindent Note that any $\xi\in \Gamma $ can be written as
$\xi=\xi^{(0)} \textbf{e}(0)\oplus \xi^{(1)} \oplus \cdots,~
\xi^{(n)}$ in the $n$-fold symmetric tensor product  $L^2(\mathbb
R_+,\mathbf k)^{\otimes n}\equiv L^2(\Sigma_n) \otimes \mathbf
k^{\otimes n}$ where $\Sigma_n$ is the $n$-simplex $\{
\underbar{t}=(t_1,t_2,\cdots, t_n):0\le t_1 <t_2\cdots <t_n
<\infty\}.$

\blema \label{qs-esti}  For any $ u\in \mathbf h,v\in \mathcal D,
0\le s\le t\le \tau$
 \be \|\sum_{j\ge 1} \int_s^t
V_{s,\lambda} (u,L_jv) a_j^\dag (d\lambda) \textbf{e}(0)\|^2\le
C_\tau (t-s)\| u\|^2 \sum_{j\ge 1} \|L_jv \|^2 \ee where $C_\tau=2
e^\tau$ \elema
\begin{proof} For  any  $\phi$ in the Fock space $\Gamma(L^2(\mathbb R_+, \mathbf k))$,
\bean && \langle \phi,\sum_{j\ge 1} \int_s^t V_{s,\lambda} (u,L_jv)
a_j^\dag (d\lambda) \textbf{e}(0)\rangle |^2\\
&&=|\langle u\otimes\phi,\{ \sum_{j\ge 1} \int_s^t V_{s,\lambda}
a_j^\dag (d\lambda) \}L_jv \otimes  \textbf{e}(0)\rangle |^2\\
&&\le \| u\otimes\phi\|^2 \|\{ \sum_{j\ge 1} \int_s^t V_{s,\lambda}
a_j^\dag (d\lambda)\} L_jv \otimes  \textbf{e}(0) \|^2. \eean By
estimate of quantum stochastic integration  (Proposition 27.1,
\cite{krp}), the above quantity is \bean  &&\le C_\tau \|
u\otimes\phi\|^2  \sum_{j\ge 1}  \int_s^t
\|V_{s,\lambda} L_jv \otimes  \textbf{e}(0) \|^2 ~d\lambda\\
&&\le C_\tau (t-s)\| u\otimes\phi\|^2  \sum_{j\ge 1} \|L_jv \|^2 .
\eean Since $\phi$  is arbitrary requirement follows.

\end{proof}
 \blema \label{lemma-vst-1} For any $ u\in \mathbf h,v\in \mathcal D
, 0\le s\le t\le \tau$  there exist constants $C_{\tau,u,v},
C_{\tau,u,v}^\prime $ given by
\[ C_{\tau,u,v}= 2 \|u\|^2[ C_\tau  \sum_{j\ge 1}
  \|L_j v\|^2 +  \tau \|G~v\|^2 ]\]
  and
  \[C_{\tau,u,v}^\prime=2 C_\tau \|u\|^2
      [ C_\tau  \sum_{i,j\ge 1}
  \|L_j L_i v\|^2 +  \tau \sum_{i\ge 1}\|G~L_i ~v\|^2 ]\]
   such that
\begin{description}

\item[(i)] $\|(V_{s,t}-1)(u,v)~~\textbf{e}(0)\|^2\le  C_{\tau,u,v}
(t-s)$

\item[(ii)] $ \|\zeta(s,t,u,v) \|^2\le C_{\tau,u,v}^\prime$ and
 $\|\varsigma(s,t,u,v)\|\le C_{\tau,u,v}  \sqrt{t-s},~ \forall ~0\le s < t\le
\tau.$

\item[(iii)]For any $\xi\in \Gamma(L^2(\mathbb R_+, \mathbf k)),$
~ $\lim_{s \rightarrow t}\langle \xi,  \zeta(s,t,u,v)\rangle =0$ and
\[\lim_{s \rightarrow t} \langle \xi, \gamma(s,t,u,v)\rangle =\sum_{j\ge 1}
\langle u,L_jv\rangle \overline{\xi^{(1)}_j} (t)= \langle \xi^{(1)}
(t), \eta(u,v)\rangle,~~\mbox{a.e.}~~ t\ge 0.\]
\end{description}
\elema

\begin{proof} {\bf (i)}  By identity  (\ref{vst-1}) and Lemma \ref{qs-esti} we have
   \bean
  && \|(V_{s,t}-1)(u,v)~~\textbf{e}(0)\|^2 \\
   &&=\|\sum_{j\ge 1} \int_s^t
  V_{s,\alpha}(u,L_j v) a_j^\dag(d\alpha)~\textbf{e}(0)+ \int_s^t
  V_{s,\alpha}(u, G v)~\textbf{e}(0) d\alpha\|^2\\
  && \le 2 \|\sum_{j\ge 1} \int_s^t
  V_{s,\alpha}(u,L_j v) d\alpha~~\textbf{e}(0)\|^2+[\int_s^t
  \|V_{s,\alpha}(u, G v)~~\textbf{e}(0)\| d\alpha]^2\\
   && \le 2 \|u\|^2[ C_\tau (t-s) \sum_{j\ge 1}
  \|L_j v\|^2 +  [(t-s)\|G~v\|]^2 ] \\
     && \le  C_{\tau,u,v} (t-s).
  \eean

\noindent {\bf(ii) }  {1.} As in the proof of  Lemma \ref{qs-esti}
we have
  \bean
  && \|\zeta(s,t,u,v)\|^2=\frac{1}{(t-s)^2}\|\sum_{j\ge 1} \int_s^t
  (V_{s,\lambda}-1)(u,L_jv) a_j^\dag (d\lambda) ~~\textbf{e}(0)\|^2\\
    && \le  \frac{\|u\|^2}{(t-s)^2}\| \sum_{j\ge 1} \int_s^t
  (V_{s,\lambda}-1)L_jv~~a_j^\dag (d\lambda) ~~\textbf{e}(0)\|^2.\eean
   Since $L_j v\in \mathcal D$  for all $j\ge 1$  by assumption {\bf E3}, by
estimate of  quantum stochastic integration  (Proposition 27.1,
\cite{krp})  the above quantity is
 \bean  &&\le  \frac{C_\tau \|u\|^2}{(t-s)^2} \sum_{j\ge 1} \int_s^t\|
  (V_{s,\lambda}-1)L_jv ~~\textbf{e}(0)\|^2 d\lambda\\
     && \le 2 \frac{C_\tau \|u\|^2}{(t-s)^2}\sum_{j\ge 1} (t-s)
      [ C_\tau (t-s) \sum_{i\ge 1}
  \|L_i L_j v\|^2 +  (t-s)^2\|G~L_j ~v\|^2 ] \\
  && \le 2 C_\tau \|u\|^2 \sum_{j\ge 1}
      [  C_\tau \sum_{i\ge 1}
  \|L_i L_j v\|^2 +  (t-s) \|G~L_j ~v\|^2 ]\\
  && \le 2 C_\tau \|u\|^2 \sum_{i\ge 1}
      [ C_\tau \sum_{j\ge 1}
  \|L_j L_i v\|^2 +  \tau \|G~L_i ~v\|^2 ]=C_{\tau,u,v}^\prime  \eean

\noindent {2.} We have
  \bean
  && \|\varsigma(s,t,u,v)\|=\frac{1}{(t-s)}\| \int_s^t
  (V_{s,\lambda}-1)
(u,Gv) d \lambda ~~\textbf{e}(0)\|\\
    && \le  \frac{1}{(t-s)} \int_s^t
  \|(V_{s,\lambda}-1)
(u,Gv)~~\textbf{e}(0)\| d\lambda.
  \eean
By part {\bf(i)} it follows that $\|\varsigma(s,t,u,v)\|^2\le
C_{\tau,u,v}  \sqrt{t-s} .$

\noindent  {\bf(iii)}
 {1.} For any
$f\in L^2(\mathbb R_+,\mathbf k)$ let us consider
 \bean
 &&\langle
\textbf{e}(f),\zeta(s,t,u,v)\rangle=\langle
\textbf{e}(f),\frac{1}{t-s} \sum_{j\ge 1}\int_s^t (V_{s,\lambda}-1)
(u,L_jv)
a_j^\dag (d\lambda) ~~ \textbf{e}(0)\rangle\\
&&=\frac{1}{t-s} \sum_{j\ge 1}\int_s^t \overline{f_j(\lambda)}
\langle \textbf{e}(f),(V_{s,\lambda}-1) (u,L_jv)  ~~
\textbf{e}(0)\rangle  d\lambda\\
&&=\frac{1}{t-s} \int_s^t G(s,\lambda) d\lambda, \eean where
$G(s,\lambda)= \sum_{j\ge 1} \overline{f_j(\lambda)} \langle
\textbf{e}(f),(V_{s,\lambda}-1) (u,L_jv)  ~~ \textbf{e}(0)\rangle. $
 Note that the complex valued function
$G(s,\lambda)$ is uniformly continuous in both the variable
$s,\lambda$ on $[0,\tau]$  and $G(t,t)=0.$ So we get
\[\lim_{s \rightarrow t} \langle
\textbf{e}(f),\zeta(s,t,u,v)\rangle=0.\] Since  $\zeta(s,t,u,v)$
uniformly bounded in $s,t$
\[\lim_{s \rightarrow t} \langle
\xi,\zeta(s,t,u,v)\rangle=0, \forall \xi \in \Gamma.\]

\noindent  {2.} We have
 \be \label{xi1}  \langle \xi,
\gamma(s,t,u,v)\rangle = \frac{1}{t-s}\sum_{j\ge 1} \langle
u,L_jv\rangle \int_s^t \overline{\xi^{(1)}_j} (\lambda )
d\lambda.\ee Since
\[|\sum_{j\ge 1} \langle
u,L_jv\rangle \overline{\xi^{(1)}_j} (t) |^2\le \|u\|^2 \sum_{j\ge
1} \|L_jv\|^2  \sum_{j\ge 1}  |\xi^{(1)}_j (t)|^2\le \sum_{j\ge 1}
\|L_jv\|^2 \|\xi^{(1)} (t)\|^2,\] the function $\sum_{j\ge 1}
\langle u,L_jv\rangle \overline{\xi^{(1)}_j} (\cdot)\in L^2$  and
hence locally integrable. Thus we get
\[\lim_{s \rightarrow t} \langle \xi,
\gamma(s,t,u,v)\rangle=\sum_{j\ge 1} \langle u,L_jv\rangle
\overline{\xi^{(1)}_j} (t)~~\mbox{a.e.}~~ t\ge 0.\]

\end{proof}
\blema \label{Mst}
 For $n\ge 1, ~ \underbar{t}\in \Sigma_n$  and $u_k\in \mathbf h,v_k\in
\mathcal D:k=1,\cdots, n, \xi \in \Gamma(L^2(\mathbb R_+,\mathbf
k))$ and $[s_k,t_k)$'s are disjoint..
\begin{description}
\item[(i)] $\lim_{\underbar{s} \rightarrow \underbar{t}}\langle \xi,
\prod_{k=1}^n M(s_k,t_k,u_k,v_k)~~\textbf{e}(0) \rangle=0,$ \\
where
$M(s_k,t_k,u_k,v_k)=\frac{(V_{s_k,t_k}-1)}{t_k-s_k}(u_k,v_k)-\langle
u_k, G~v_k \rangle-\gamma(s_k,t_k,u_k,v_k)$  and $\lim_{\underbar{s}
\rightarrow \underbar{t}}$  means  $s_k\rightarrow t_k$  for each
$k.$

\item[(ii)]
$\lim_{\underbar{s} \rightarrow \underbar{t}}\langle \xi,
\otimes_{k=1}^n \gamma(s_k,t_k,u_k,v_k) \rangle=\langle
\xi^{(n)}(t_1,\cdots, t_n), \eta(u_1,v_1) \otimes\cdots \otimes
\eta(u_n,v_n) \rangle.$
\end{description}
 \elema
  \begin{proof}
 {\bf(i)} First note that
$M(s,t,u,v)\textbf{e}(0)=\zeta(s,t,u,v)+~ \varsigma(s,t,u,v).$
 So by the above observations
 $\{M(s,t,u,v) \textbf{e}(0)\}$ is  uniformly bounded in $s,t$
 and\\
$\lim_{s \rightarrow t} \langle \textbf{e}(f)
,M(s,t,u,v)\textbf{e}(0)\rangle=0, \forall f\in L^2(\mathbb
R_+,\mathbf k).$ Since the intervals $[s_k,t_k)$'s are disjoint for
different $k$'s,
\[\langle \textbf{e}(f), \prod_{k=1}^n M(s_k,t_k,u_k,v_k)~~\textbf{e}(0)
\rangle= \prod_{k=1}^n \langle \textbf{e}(f_{[s_k,t_k)}),
M(s_k,t_k,u_k,v_k)~~\textbf{e}(0) \rangle\] and thus
$\lim_{\underbar{s} \rightarrow \underbar{t}}\langle \textbf{e}(f),
\prod_{k=1}^n M(s_k,t_k,u_k,v_k)~~\textbf{e}(0) \rangle=0.$ By Lemma
\ref{lemma-vst-1}, the vector  $\prod_{k=1}^n
M(s_k,t_k,u_k,v_k)~~\textbf{e}(0)$ is uniformly bounded  in
$s_k,t_k$
and hence  convergence hold if we replace  $\textbf{e}(f)$  by   any  vector $\xi$ in  the Fock Space.\\

\noindent {\bf(ii) } It can be proved similarly as in part {\bf(iii)
} of the  previous Lemma.
\end{proof}

\blema Let $\xi\in \Gamma$ be such that  \be \label{ortho} \langle
\xi, \zeta\rangle=0,~\forall \zeta \in \mathcal S^\prime, \ee Then
\begin{description}
\item[(i)]$\xi^{(0)}=0$   and  $\xi^{(1)}(t)=0$
~\mbox{for  a.e.} $t\in[0,\tau].$
\item[(ii)]For any $n\ge 0,~ \xi^{(n)}( \underbar{t})=0$
~\mbox{for  a.e.} $\underbar{t}\in \Sigma_n: t_i\le \tau.$
\item[(iii)]The set  $\mathcal S^\prime$ is total in the  Fock space $\Gamma.$
\end{description}
\elema

\begin{proof}
{\bf(i)} For any $s\ge 0,~V_{s,s}=1_{\mathbf h\otimes \Gamma}$ so in
particular (\ref{ortho}) gives, for any  $u\in \mathbf h, v \in
\mathcal D$
\[0=\langle \xi, V_{s,s}(u,v)  \textbf{e}(0)\rangle=\langle u,v\rangle
\overline{\xi^{(0)}}\] and hence $\xi^{(0)}=0.$\\

\noindent By (\ref{ortho}),~ $\langle \xi, [V_{s,t}-1] (u,v)
\textbf{e}(0)\rangle=0$ for any  $0\le s < t\le \tau<\infty , u\in
\mathbf h, v\in \mathcal D.$ Hence for any  $ u\in \mathbf h, v\in
\mathcal D$  by
 Lemma \ref{vst-1} (iii) we have  \bean &&0=\lim_{s \rightarrow t}
\frac{1}{t-s}\langle \xi, [V_{s,t}-1] (u,v) \textbf{e}(0)\rangle \\
&&=\sum_{j\ge 1}\langle u,L_jv \rangle \overline{\xi^{(1)}_j(t)}
=\sum_{j\ge 1}\eta_j(u,v) \overline{\xi_j^{(1)}(t)}=\langle
\xi^{(1)}(t), \eta(u,v) \rangle \eean for almost all $t\in
[0,\tau].$
 Since $\{\eta(u,v):u\in \mathbf h, v\in \mathcal D\}$ is total in
$\mathbf k$ it follows that $\xi^{(1)}(t)=0$ for almost all  $ t\le
\tau.$

\noindent {\bf(ii)} We prove this by induction. The result is
already proved for $n=0,1.$ For $n\ge 2,$ assume as induction
hypothesis  that for all $m\le n-1,$ $\xi^{(m)}( \underbar{t})=0,$
~\mbox{for a.e.} $\underbar{t}\in \Sigma_m: t_k\le \tau,
 k=1,2,\cdots, m.$ We now   show that $\xi^{(n)}( \underbar{t})=0,$
~\mbox{for
a.e.} $\underbar{t}\in \Sigma_n: t_k\le \tau.$\\

 \noindent Let $0\le s_1< t_1\le s_2<t_2<\ldots < s_n <t_n\le \tau$
 and $u_k\in \mathbf h, v_k\in \mathcal D:k=1,2\cdots,n.$
By (\ref{ortho})  and  part ({\bf i})  we have
\[\langle \xi,
 \prod_{k=1}^n \frac{(V_{s_k,t_k}-1)}{t_k-s_k}(u_k,v_k)~~\textbf{e}(0) \rangle =0.\]
Thus
 \bea \label{nvst-1}
 && 0=\lim_{\underbar{s}
\rightarrow \underbar{t}}\langle \xi,
 \prod_{k=1}^n \frac{(V_{s_k,t_k}-1)}{t_k-s_k}(u_k,v_k)~~\textbf{e}(0) \rangle \\
 && =\lim_{\underbar{s}
\rightarrow \underbar{t}}\langle \xi, \prod_{k=1}^n \{
M(s_k,t_k,u_k,v_k)+\langle u_k, G~v_k \rangle+
 \gamma(s_k,t_k,u_k,v_k) \}~~\textbf{e}(0) \rangle. {\nonumber}
\eea Let $P,Q,R$   and $ P^\prime, R^\prime $ be two sets of
disjoint  partitions of $ \{1,2,\cdots , n\}$ such that $Q$ and $ R$
are non empty. We write $|S|$ for the cardinality  of set  $S.$ Then
by Lemma \ref{Mst} (ii) the  right hand side of (\ref{nvst-1}) is
equal to
 \bean
 && \sum_{P^\prime
,R^\prime }  \langle
\xi^{(|R^\prime|)}(t_{r_1^\prime},\cdots,t_{r_{|R^\prime|}^\prime
}), \otimes_{k \in R^\prime }~\eta(u_k,v_k) \rangle ~\prod_{k\in
P^\prime
} \langle u_k, G~v_k \rangle\\
 &&~~~ +  \lim_{\underbar{s}
\rightarrow \underbar{t}}\sum_{P,Q,R}\langle \xi,  \prod_{k\in P}
\langle u_k, G~v_k \rangle~\prod_{k\in Q}\{ M
(s_k,t_k,u_k,v_k)\}\prod_{k \in R} \{\gamma(s_k,t_k,u_k,v_k)
\}~~\textbf{e}(0) \rangle.
 \eean
Thus by the induction hypothesis,
 \bea \label{lastterm}&& 0= \langle  \xi^{(n)}(t_1,t_2,\cdots, t_n),
\eta(u_1,v_1)
\otimes\cdots \otimes  \eta(u_n,v_n) \rangle\\
 &&~~~ +  \lim_{\underbar{s}
\rightarrow \underbar{t}}\sum_{P,Q,R}\langle \xi,  \prod_{k\in P}
\langle u_k, G~v_k \rangle~\prod_{k\in Q}\{ M
(s_k,t_k,u_k,v_k)\}\prod_{k \in R} \{\gamma(s_k,t_k,u_k,v_k)
\}~~\textbf{e}(0) \rangle. {\nonumber}
 \eea
We claim that the second  term in  (\ref{lastterm}) vanishes. To
prove the claim, it is enough to  show  that for any two non empty
disjoint subsets $Q\equiv\{q_1,q_2,\cdots,
q_{|Q|}\},R\equiv\{r_1,r_2, \cdots, r_{|R|}\}$  of   $ \{1,2,\cdots,
n\},$
 \be \label{QR00}
 \lim_{\underbar{s}
\rightarrow \underbar{t}} \langle \xi, \prod_{q\in Q}\{ M
(s_q,t_q,u_q,v_q)\}\prod_{r \in R} \{\gamma(s_r,t_r,u_r,v_r)
\}~~\textbf{e}(0) \rangle=0.\ee Writing $ \psi $ for the vector
$\prod_{q\in Q}\{ M (s_q,t_q,u_q,v_q) \}\textbf{e}(0),$ we have
 \bea \label{QR11}
&& \langle \xi, \prod_{q\in Q}\{ M (s_q,t_q,u_q,v_q)\}\prod_{r \in
R} \{\gamma(s_r,t_r,u_r,v_r)
\}~~\textbf{e}(0) \rangle    {\nonumber}\\
&&=  \langle \xi, \psi \otimes  \otimes_{r \in
R} \frac{1_{[s_r,t_r]}~\eta(u_r,v_r)}{t_r-s_r} \rangle   {\nonumber}\\
&&=\langle \xi,\psi \otimes \otimes_{r \in R}
\frac{1_{[s_r,t_r]}~\eta(u_r,v_r)}{t_r-s_r}  \rangle {\nonumber}\\
&&=  \sum_{l\ge |R|} \langle \xi^{(l)},\psi^{(l-|R|)} \otimes
\otimes_{r \in R} \frac{1_{[s_r,t_r]}~\eta(u_r,v_r)}{t_r-s_r}
\rangle {\nonumber}\\
 &&=\langle \sum_{l\ge |R|} \langle \langle
\xi^{(l)},\psi^{(l-|R|)} \rangle \rangle, \otimes_{r \in R}
\frac{1_{[s_r,t_r]}~\eta(u_r,v_r)}{t_r-s_r}  \rangle. \eea
 Here  $ \langle \langle
\psi^{(l-|R|)}, \xi^{(l)} \rangle \rangle \in L^2(\mathbb
R_+,\mathbf k)^{\otimes |R|}$ is  defined  as in (\ref{partinn}) by
 \bea
\label{<<>>} && \langle ~~\langle \langle \psi^{(l-|R|)}, \xi^{(l)}
\rangle \rangle, \rho^{(|R|)}~~ \rangle =\langle
\xi^{(l)},\psi^{(l-|R|)}
\otimes  \rho^{(|R|)} \rangle\\
&&=\int_{\Sigma_l}\langle \xi^{(l)} (x_1,x_2,\cdots,
x_l),  {\nonumber}\\
&&~~~~~~~~~~~~~~~~~~~~~~\psi^{(l-|R|)} (x_1, x_2,\cdots,
x_{l-|R|})\otimes \rho^{(|R|)} (x_{l-|R|+1},\cdots, x_l)
\rangle_{\mathbf k^{\otimes l}}~~ dx {\nonumber} \eea for any
$\rho^{(|R|)}\in L^2(\mathbb R_+,\mathbf k)^{\otimes |R|}.$\\
 By
Lemma \ref{Mst} (i),
 \be
  \lim_{s_q \rightarrow t_q}
  \langle \xi, \prod_{q\in Q}\{ M
(s_q,t_q,u_q,v_q)\}\prod_{r \in R} \{\gamma(s_r,t_r,u_r,v_r)
\}~~\textbf{e}(0) \rangle=0.\ee  However, we need to prove
(\ref{QR00}) where the limit $\underbar{s} \rightarrow \underbar{t}$
has to be in arbitrary order. On the other hand, by (\ref{QR11}) and
(\ref{<<>>}) we get
 \bea \label{QR33} && \lim_{s_q \rightarrow t_q} \lim_{s_r \rightarrow t_r}\langle \xi,
\prod_{q\in Q}\{ M (s_q,t_q,u_q,v_q)\}\prod_{r \in R}
\{\gamma(s_r,t_r,u_r,v_r) \}~~\textbf{e}(0) \rangle
{\nonumber}\\
&&= \lim_{s_q \rightarrow t_q} \lim_{s_r \rightarrow t_r} \langle
\sum_{l\ge |R|} \langle \langle \psi^{(l-|R|)},  \xi^{(l)} \rangle
\rangle, \otimes_{r
\in R} \frac{1_{[s_r,t_r]}~\eta(u_r,v_r)}{t_r-s_r}  \rangle{\nonumber}\\
&&= \lim_{s_q \rightarrow t_q} \lim_{s_r \rightarrow t_r} \langle
\int_{\Sigma_{|R|}}\langle [\sum_{l\ge |R|}  \langle \langle
\psi^{(l-|R|)} , \xi^{(l)} \rangle \rangle ](x_1,x_2,\cdots, x_{|R|}), {\nonumber}\\
&&~~~~~~~~~~~~~~~~~~~~~~~ \otimes_{r \in R} \frac{1_{[s_r,t_r]}
(x_r)~~\eta(u_r,v_r)}{t_r-s_r} \rangle dx{\nonumber}\\
&&= \lim_{s_q \rightarrow t_q}   \langle \sum_{l\ge |R|} \langle
\langle \psi^{(l-|R|)} ,  \xi^{(l)}\rangle \rangle(t_{r_1}, \cdots,
t_{r_{|R|}}), \otimes_{r \in R}~\eta(u_r,v_r) \rangle,\eea for
almost all $\underbar{t}\in \Sigma_{|R|}.$ We fix $\underbar{t}\in
\Sigma_{|R|}$ and define families of vectors  $\widetilde{\xi}^{(l)}
:l\ge 0$  in $ L^2(\mathbb R_+, \mathbf k) ^{\otimes l}$ by
 \bean && \widetilde{\xi}^{(0)} = \langle \xi^{(|R|)}(t_{r_1}, \cdots,
t_{r_{|R|}}), \otimes_{r \in
R}~\eta(u_r,v_r) \rangle \in \mathbb C\\
&& \widetilde{\xi}^{(l)} (x_1,x_2,\cdots , x_l) = \langle\langle
\otimes_{r \in R}~\eta(u_r,v_r) ,\xi^{(|R|+l)}(x_1,\cdots ,x_l,
t_{r_1}, \cdots, t_{r_{|R|}})\rangle \rangle, \eean which defines a
Fock space vector $\widetilde{\xi}.$ Therefore, from (\ref{QR33}),
we get  that \bean && \lim_{s_q \rightarrow t_q} \lim_{s_r
\rightarrow t_r} \langle \xi, \prod_{q\in Q}\{ M
(s_q,t_q,u_q,v_q)\}\prod_{r \in R} \{\gamma(s_r,t_r,u_r,v_r)
\}~~\textbf{e}(0) \rangle = \lim_{s_q
\rightarrow t_q}   \langle \widetilde{\xi}~,~\psi \rangle\\
&&= \lim_{s_q \rightarrow t_q}   \langle \widetilde{\xi}~,~
[\prod_{q\in Q} M (s_q,t_q,u_q,v_q)] ~\textbf{e}(0)\rangle,\eean
which is equal to $0$ by Lemma \ref{Mst}  (i).
 Thus from (\ref{lastterm})  we get that  \[\langle
\xi^{(n)}(t_1,t_2,\cdots, t_n), \eta(u_1,v_1) \otimes\cdots \otimes
\eta(u_n,v_n) \rangle=0.\] Since $\{ \eta(u,v):u\in \mathbf h, v\in
\mathcal D\}$ is total in $\mathbf k,$ it follows that
$\xi^{(n)}(t_1,t_2,\cdots, t_n)=0$ for almost every
$(t_1,t_2,\cdots, t_n)\in \Sigma_n:t_k\le \tau.$
\end{proof}

\noindent {\bf(iii)} Since $\tau\ge 0 $  is arbitrary
$\xi^{(n)}=0\in L^2(\mathbb R_+, \mathbf k)^{\otimes n}:n\ge 0$  and
hence $\xi=0.$ Which proves the  totality  of $\mathcal
S^\prime\subseteq  \Gamma.$

\subsection{Unitary Equivalence} Here  we shall prove the part  (ii) of the {\bf
Theorem \ref{mainthm} } that the unitary evolution $\{U_{t}\}$ on
$\mathbf h \otimes \mathcal H $ is unitarily equivalent to the
unitary solution  $\{V_{t}\}$ of HP equation (\ref{hpeqn,st11}). To
prove this we need the following two results. Let us recall that the
subset $ \mathcal S=\{ \xi=
U_{\underbar{s},\underbar{t}}(\underbar{u},\underbar{v})
\Omega:=U_{s_1,t_1}(u_1,v_1)\cdots U_{s_n,t_n}(u_n,v_n)\Omega :
\underbar{s}=(s_1,s_2, \cdots, s_n), \underbar{t}=(t_1,t_2, \cdots,
t_n)$ $:~ 0 \le s_1\le t_1\le \cdots \le s_n\le  t_n< \infty,n\ge 1,
\underbar{u}=\otimes_{i=1}^n u_i\in \mathbf h^{\otimes
n},\underbar{v}=\otimes_{i=1}^n v_i\in \mathcal D^{\otimes n}\}$ is
total in $\mathcal H$
 and the subset\\
$ \mathcal S^\prime:=\{ \zeta=
V_{\underbar{s},\underbar{t}}(\underbar{u},\underbar{v})
\textbf{e}(0):=V_{s_1,t_1}(u_1,v_1)\cdots
V_{s_n,t_n}(u_n,v_n)\textbf{e}(0): \underbar{u}=\otimes_{i=1}^n
u_i\in \mathbf h^{\otimes n},\underbar{v}=\otimes_{i=1}^n v_i\in
\mathcal D^{\otimes n},\underbar{s}=(s_1,s_2, \cdots, s_n),
\underbar{t}=(t_1,t_2, \cdots, t_n)\}$ is total in $\Gamma.$ \blema
Let  $ U_{\underbar{s},\underbar{t}}(\underbar{u},\underbar{v})
\Omega, ~~
U_{\underbar{s}^\prime,\underbar{t}^\prime}(\underbar{p},\underbar{w})
\Omega \in \mathcal S.$\\
  Then there exist an integer $m\ge 1,~\underbar{a}=(a_1,a_2,
\cdots, a_m), \underbar{b}=(b_1,b_2, \cdots, b_m):~ 0 \le a_1\le
b_1\le  \cdots\le a_m\le  b_m< \infty,$  partition  $R_1\cup R_2
\cup R_3=\{1,\cdots, m\}$ with $|R_i|=m_i,$  family of vectors
$x_{k_l} , g_{k_i} \in \mathbf h $ and $y_{k_l}, h_{k_i}  \in
\mathcal D :l\in R_1 \cup R_2,~  i\in R_2 \cup R_3$ such that \be
 U_{\underbar{s},\underbar{t}}(\underbar{u},\underbar{v})
=\sum_{\underbar{k}} \prod_{l\in R_1 \cup R_2} U_{a_l,b_l}
(x_{k_l},y_{k_l}) \ee

\be
 U_{\underbar{s}^\prime,\underbar{t}^\prime}(\underbar{p},\underbar{w})
=\sum_{\underbar{k}} \prod_{l\in R_2 \cup R_3} U_{a_l,b_l}
(g_{k_l},h_{k_l}). \ee \elema
\begin{proof}
 It follows  from the evolution hypothesis of  the  family  of
unitary  operators  $\{U_{s,t}\}$  as for  $r\in [s,t]$  and
orhonormal basis  $\{f_j\}\subseteq \mathcal D$ of $\mathbf h$ we
can write $U_{s,t}(u,v)= \sum_{j\ge 1} U_{s,r} (u,
f_j)U_{r,t}(f_j,v).$
\end{proof}
\brmrk \label{V} Since the family of unitaries  $\{  V_{s,t}\}$ on
$\mathbf h \otimes \Gamma$ enjoy all the properties  satisfy by
family of unitaries  $\{  U_{s,t}\}$ on $\mathbf h \otimes \mathcal
H$ the above Lemma also hold  if we replace $U_{s,t}$ by $V_{s,t}.$
\ermrk
 \blema  \label{uvinner} For   $
U_{\underbar{s},\underbar{t}}(\underbar{u},\underbar{v}) \Omega, ~~
U_{\underbar{s}^\prime,\underbar{t}^\prime}(\underbar{p},\underbar{w})
\Omega \in \mathcal S.$\\
\be \label{uvinn} \langle
U_{\underbar{s},\underbar{t}}(\underbar{u},\underbar{v}) \Omega,
U_{\underbar{s}^\prime,\underbar{t}^\prime}(\underbar{p},\underbar{w})\Omega
\rangle=\langle
V_{\underbar{s},\underbar{t}}(\underbar{u},\underbar{v})
\textbf{e}(0),
V_{\underbar{s}^\prime,\underbar{t}^\prime}(\underbar{p},\underbar{w})
\textbf{e}(0) \rangle.
 \ee
  \elema
\begin{proof}
We have by previous Lemma  and  assumption {\bf A}
 \bean &&\langle
U_{\underbar{s},\underbar{t}}(\underbar{u},\underbar{v}) \Omega,
U_{\underbar{s}^\prime,\underbar{t}^\prime}(\underbar{p},\underbar{w})\Omega
\rangle\\
&&=\sum_{\underbar{k}} \prod_{l\in R_1} \langle U_{b_l-a_l}
(x_{k_l},y_{k_l})\Omega ,\Omega\rangle \prod_{l\in R_2} \langle
U_{b_l-a_l} (x_{k_l},y_{k_l})\Omega,
 U_{b_l-a_l} (g_{k_l},h_{k_l}) \Omega \rangle \\
 &&~~~~~~~~~~
 \prod_{l\in R_3} \langle
\Omega , U_{b_l-a_l} (g_{k_l},h_{k_l}) \Omega \rangle \\
&&=\sum_{\underbar{k}} \prod_{l\in R_1} \langle T_{b_l-a_l}
y_{k_l},x_{k_l}\rangle \prod_{l\in R_2} \langle g_{k_l} ,
 Z_{b_l-a_l} (| h_{k_l} >< y_{k_l} |)~x_{k_l} \rangle
 \prod_{l\in R_3} \langle
g_{k_l}  , T_{b_l-a_l} h_{k_l}\rangle \\
&&=\sum_{\underbar{k}} \prod_{l\in R_1} \langle V_{b_l-a_l}
(x_{k_l},y_{k_l})\textbf{e}(0), ~\textbf{e}(0)\rangle \prod_{l\in
R_2} \langle V_{b_l-a_l} (x_{k_l},y_{k_l})\textbf{e}(0),
 V_{b_l-a_l} (g_{k_l},h_{k_l})\textbf{e}(0) \rangle\\
 &&~~~~~~~~~~
 \prod_{l\in R_3} \langle
\textbf{e}(0) , V_{b_l-a_l} (g_{k_l},h_{k_l}) \textbf{e}(0)\rangle.
\eean Now  by Remark (\ref{V}), the above quantity is equal to
  $\langle
V_{\underbar{s},\underbar{t}}(\underbar{u},\underbar{v})
\textbf{e}(0),
V_{\underbar{s}^\prime,\underbar{t}^\prime}(\underbar{p},\underbar{w})
\textbf{e}(0) \rangle.$
\end{proof}

\noindent {\bf Proof of the  part (ii) of Theorem \ref{mainthm}~:} \\
We need to construct a unitary operator  $\widetilde{\Xi}:\mathbf h
\otimes \mathcal H \rightarrow \mathbf h \otimes \Gamma$  such that
\be \label{U=V} U_t =\widetilde{\Xi}^*~ V_t~
\widetilde{\Xi},~\forall~~ t\ge 0.\ee
 Let
us define a map $\Xi: \mathcal H \rightarrow \Gamma$ by setting, for
any $\xi= U_{\underbar{s},\underbar{t}}(\underbar{u},\underbar{v})
\Omega \in \mathcal S,$~~  $\Xi
\xi:=V_{\underbar{s},\underbar{t}}(\underbar{u},\underbar{v})
\textbf{e}(0)\in \mathcal S^\prime$ and then extending linearly. So
by definition and totality of $\mathcal S^\prime,$ range of $\Xi$ is
dense in $\Gamma.$
 To see that $\Xi$ is   a unitary  operator from $\mathcal H$ to
 $\Gamma$ it is enough to note from Lemma \ref{uvinner}  that
\be \langle   \Xi \xi, \Xi \xi^\prime  \rangle=\langle   \xi,
\xi^\prime  \rangle, ~ \forall~~ \xi,\xi^\prime \in \mathcal S.\ee
For the conclusion it is suffices to set $\widetilde{\Xi}=1_{\mathbf
h} \otimes \Xi .$ \qed

\brmrk The  assumption {\bf C}  is  ruling  out the presence  of
conservation  (Poisson )terms  in the  associated  HP equation  as
the representation $\pi,$ we obtained,  is  trivial (see Remark
\ref{trivial-rep}). Without this   assumption {\bf C}, the problem
is not yet settled. In the absence of  assumption {\bf C} the
representation $\pi$ shall be  non trivial which in general will
give  rise  to a unitary (different from identity) operator $W$ on
$\mathbf h \otimes \mathbf k$ and associated HP equation
(\ref{hpeqn,st11}) will contain conservation terms with coefficients
$\{L_\nu^\mu\}$  described as in  (\ref{Lmunu}).

\ermrk

\brmrk  The Hypothesis  {\bf E2}, i.e. there exists
  $ \mathcal D, $  core for $G$  such that $ \mathcal D\subseteq \mathcal D(L_j^*)$ for every  $j\ge
  1,$  is a  strong assumption. But this is  necessary  one  in
  order that  quantum stochastic  differential  equation for   $V_{t}$
   makes  sense. Only way one can do away with this assumption is to
   abandon the quantum stochastic differential equation   for $V_{t}$  and just deal with $V_t$ as a left
   cocycle  described by the associated four semigroups \cite{lw2}.
   This programme is not yet complete.

\ermrk

\brmrk  The Hypothesis  {\bf E3}, i.e.  for any $v\in \mathcal D,
\sum_{j\ge 1} \| G L_j v\|^2  <\infty.$  This  holds trivially  when
$[ G ,L_j]=0.$  Condition $[ G ,L_j]=0,$   in particular holds for
classical Brownian motion on $\mathbb R^n$  and  for Casimir
operator  $G$ on Lie algebra of a locally compact Lie  group
$\mathcal G$ with $L_j=X_j$ represented  on the Hilbert space
$\mathbf h= L^2(\mathcal G),$ where $\{X_j\}_{j=1}^n $ a basis for
the Lie algebra. The commutator $[ G ,L_j]$  also vanish  in case of
Quantum Brownian motion on non-commutative Torus, Quantum Heisenberg
manifold and Quantum Plane \cite{gs} . \ermrk

\end{document}